\documentclass[12pt]{article}
\usepackage{newinutile}
\usepackage{amsmath,amssymb,esint}
\usepackage{epsfig,graphics,color,calc,graphicx}
\usepackage{siunitx}
\usepackage{mathrsfs}

\newcommand{\rr}[1]{{\normalfont\textrm{#1}}}

\newcommand{\bb}[1]{{\mathbb{#1}}}

\begin{document}
\title{Driven particle flux through a membrane: Two-scale asymptotics of a diffusion
equation with polynomial drift}

\author{Emilio N.M.\ Cirillo}
\email{emilio.cirillo@uniroma1.it}
\affiliation{Dipartimento di Scienze di Base e Applicate per
             l'Ingegneria, Sapienza Universit\`a di Roma,
             Italy.}

\author{Ida de Bonis}
\email{i.debonis@unifortunato.eu}
\affiliation{Universit\`a degli Studi ``Giustino Fortunato",
             Benevento, Italy.}

\author{Adrian Muntean}
\email{adrian.muntean@kau.se}
\affiliation{Department of Mathematics and Computer Science,
Karlstad University, Sweden.}

\author{Omar Richardson}
\email{omar.richardson@kau.se}
\affiliation{Department of Mathematics and Computer Science,
Karlstad University, Sweden.}


\begin{abstract}

Diffusion of particles through an heterogenous obstacle 
line is modeled as a two-dimensional diffusion problem with a one--directional nonlinear convective drift and is examined using two-scale asymptotic analysis. At the scale where the structure of heterogeneities is observable the obstacle line has an inherent thickness.  Assuming the heterogeneity to be made of an array of periodically arranged microstructures (e.g. impenetrable solid rectangles),  two scaling regimes are identified: the characteristic size of the microstructure is either significantly smaller than the thickness of the obstacle line or it is of the same order of magnitude.
We  scale up the convection-diffusion model  
and compute the effective diffusion and drift tensorial coefficients for both scaling regimes. The upscaling procedure combines ideas of two-scale asymptotics homogenization with dimension reduction arguments. Consequences of these results for the construction of more general transmission boundary conditions are discussed. We numerically illustrate the behavior of the upscaled membrane in the finite thickness regime and apply it to describe the transport of {\rm CO}$_2$ through paperboard. 
\end{abstract}


\preprint{Appunti: \today}


\maketitle


\section{Introduction}
\label{s:intro}
\par\noindent

The study of the physics of interfaces has known a 
great impulse in the last decades \cite{NPW}, mainly motivated by the 
study of surfaces separating two different phases. 
Interface fluctuations, controlled by surface tension, have been studied 
with the methods of statistical mechanics, in particular those borrowed from
the theory of equilibrium critical phenomena. 
Membrane--like interfaces, namely, surfaces made of a different 
kind of molecules with respect to those forming the medium, 
do not need to separate regions of space filled with different phases, but 
they exhibit wide fluctuations, too, due to the smallness of their 
surface tension. In particular, depending on the temperature, 
they can undergo a phase transition between a flat and a
crumpled phase \cite{CGP}.

In this paper we investigate flat static (not fluctuating) membranes 
separating two regions of space and crossed by a fluid. This is the typical 
setup one is interested in when studying membrane filtration. Traditionally,  membrane filtration is one of the most common methods for purifying fluids; see e.g. \cite{Filtration} and references cited therein.
Furthermore, recent advances in conductive and mass transport through a composite medium have led to increased interest in the process of mixed-matrix membrane separation. In such cases, small particles of a microporous material, identified
as a filler, are dispersed in a dense nonporous polymer material, identified as a matrix, and then processed into a thin
composite layer, identified as a membrane. The objective is that the filler, chosen for its high adsorption
affinity or transport rate for a molecular species of interest, improves the efficacy of the matrix in membrane-mediated
separation \cite{Poz}. Depending on pore sizes and level of microscopic activity, one also encounters the so-called enhanced matrix diffusion \cite{Sato}.

Our main motivation is to develop  multiscale mathematical modelling 
strategies of transport processes that can describe, 
over several space scales, how internal structural features of the filler and local defects 
affect the permeability of the material, perceived  as a thin long permeable 
membrane. As concrete applications we have in mind the transport of {\rm O}$_2$ and  {\rm CO}$_2$  molecules through packaging materials (paperboard)  as well as the dynamics of human crowds  through barrier-like heterogeneous environments (active particles walking inside geometries with obstacles). 

We study the diffusion of particles through such a thin heterogeneous membrane 
under a one--directional nonlinear drift. Using the mean--field equation 
\begin{equation}
\label{int000}
\frac{\partial u}{\partial t}
-d_1\frac{\partial^2u}{\partial x_1^2}-d_2\frac{\partial^2u}{\partial x_2^2}
=-b\frac{\partial}{\partial x_1}[u(1-u)]+f(x),
\end{equation}
with $b>0$,
which is found in the hydrodynamic limit of the 
two--dimensional random walk with simple exclusion and drift along the $x_1$-direction
(for details, 
see \cite{CKMSSphysicaA2016}), we study the possibility to upscale the system 
and to compute the effective transport coefficients accounting for the 
presence of the membrane, adding analytic results to our simulation 
study \cite{CKMSpre2016}.

In \cite{CKMSSphysicaA2016, CKMSpre2016} the same problem is addressed in a 
microscopic setup. A lattice model, known as \emph{simple exclusion model},
is considered on a two--dimensional strip of $\mathbb{Z}^2$. 
There, particles move randomly on the strip with the constraint that 
at most one particle at a time can occupy the sites of the lattice.  
Particles move choosing at random one of the four neighboring sites and 
a \emph{drift} is introduced in the dynamics so that 
one of the four direction is possibly more probable. 
This model is a generalization of the celebrated TASEP (total 
asymmetric simple exclusion model) which is a one--dimensional 
simple exclusion model in which particles move to the right at random 
times \cite{DJLS}. 

In this framework, the equation \eqref{int000} is derived in 
the macroscopic diffusive limit, i.e., when the space and the 
drift are rescaled with a small parameter and, 
correspondingly, the time is rescaled with the square of the same 
parameter. In \cite{CKMSSphysicaA2016} we have reported a 
useful heuristic derivation of this equation which, in the 
one--dimensional case, was rigorously proven in \cite{dMPS}
(see, also, \cite{KL} for an account of the more recent techniques 
developed in the framework of hydrodynamic limit theory). 
In particular, this heuristic computation shows that the 
two diffusion coefficients can be different as a consequence of the fact 
that at the microscopic level the probability of a particle to move 
horizontally or vertically can differ. Moreover, and this is 
much more important in our context, the peculiar structure 
of the transport term on the right hand side is related to 
the probability of a particle performing a move, which the simple exclusion might prevent.
Consequently, the factor $u$ comes from the probability to find a 
particle at a given site and the factor $1-u$ accounts for 
the probability that the site where the particle tries to move to is 
indeed empty. Thus, we can summarize this discussion saying that 
the peculiar form of the right hand side of equation \eqref{int000} 
is, at the microscopic level, connected to the hard--core 
repulsion of the molecules. 

We stress that 
the model we have in mind is \eqref{int000}, but 
the techniques that will 
be developed in this article will apply to a much more general transport 
term obtained by substituting $u(1-u)$ with a general 
polynomial in terms of $u$. 


For a special scaling regime, we perform a 
simultaneous homogenization asymptotics and dimension reduction, allowing us 
not only to replace the heterogeneous membrane by an homogeneous obstacle line, 
but also to provide the effective transmission conditions needed to complete 
the upscaled model equations. The heterogeneities we account for in this 
context are assumed to be arranged periodically, but the same methodology 
can be adapted to cover also the locally periodic case. Additionally, we investigate also the effect of diffusion correlations and cross-diffusion (diagonal vs. full diffusion tensors) on the structure of the upscaled equations. We observe that in the case of the infinitely thin upscaled membrane the structure of the limit equations is unchanged, while in the case of the finite-length upscaled membrane the presence of the off-diagonal terms does not permit the use of closed form representations of oscillations in terms of cell functions. Furthermore, it is worth mentioning that the clogging of the membrane cannot be achieved with our model. Local clogging can eventually be reached by allowing the boundaries of the microstructures to evolve freely. As working techniques, 
we employ scaling arguments as well as two-scale homogenization asymptotic 
expansions to guess the structure of the model equations and the corresponding 
effective transport coefficients. As a long term plan, we would like to see whether infinitely-thin periodic membrane models can be used to give insight in the nonlinear structure of localized singularities arising in reaction terms connected to quenching structures; see for instance the settings from \cite{Ida1} and \cite{Ida2}. The question here is what a microscopic membrane would model look like so that it gives rise to production terms by reaction of the form $\eta(r,s)=k\frac{r}{s^\gamma}$ in a certain asymptotic regime, where $k>0, 0<\gamma\leq 1$ for $r,s\in [0,\infty]$ for coupled systems of semi-linear reaction-diffusion equations (cf. \cite{deBM}).

The research presented in this article pursues a formal asymptotics route; it follows the thread of the original mathematical analysis work by M.\ Neuss-Radu and 
W.\ J\"ager in \cite{membrane} by adding to the discussion the presence of 
nonlinear transport terms and is remotely related to our work on filtration 
combustion through heterogeneous thin layers; compare \cite{Ekeoma}. Recent follow-up (mathematical analysis) works of \cite{membrane} are reported in \cite{Brizzi,Gahn-Vietnam,Gahn-DCDS} (where the authors apply the concept of two-scale boundary layer convergence to the corresponding setting).  Strongly connected scenarios to the transport-through-membranes problem are the theoretical estimation of the effective interfacial resistance of regular rough surfaces (cf- \cite{Donato}, e.g.) and the upscaling of reaction, diffusion, and flow processes in porous media with thin fissures (cf. \cite{Amaziane,Zhao}, e.g.). 

What makes our study peculiar and innovative is the combination of the 
heterogeneous structure of the space region where particles move 
and the presence of the transport term on the right-hand side in 
the evolution equation \eqref{int000}.
Indeed, our results extend to a more general model 
assuming the transport term 
to be the $x_1$-derivative of a polynomial of the field $u$ with a
finite arbitrary large degree.  The main finding of this study can be summarized as follows: 
\begin{itemize}
\item We deduced the structure of the formal asymptotic expansions which are behind the concept of two-scale boundary layer convergence from \cite{membrane}; this structure can be further employed to construct corrector estimates to justify the upscaling and to provide convergence rates.
\item We derived the structure of the upscaled transmission conditions across the obstacle line with the corresponding  jumps in both transport fluxes and concentrations expressed in terms of the (local) physics taking place inside the microstructures (heterogeneities) of the membrane.  
\item Using finite element approximations of our model equations implemented in FEniCS (\cite{Logg}), we numerically illustrate the behavior of the upscaled membrane in the finite thickness regime. We simulate the basic membrane scenario using a reference parameter set corresponding to the penetration of gaseous {\rm CO}$_2$ through a porous paper sheet.    This gives confidence that our model equations and their implementation can be  used in practical applications and,  in principle, can be  extended to cover more complex membrane microstructures than the locally periodic regime. 
\end{itemize}

The article is organized as follows: In Section \ref{s:modello} we present the equations of our mean-field model as well as the membrane geometry. After a suitable scaling, we point out two relevant asymptotic regimes in terms of a small parameter $\varepsilon$ which incorporates the periodicity and selected size effects of the internal structure of the membrane. Section \ref{s:duescale} contains the derivation of the finite thickness upscaled membrane model, while in Section \ref{s:infinitesima} we consider the more delicate case of the upscaling of the infinitely-thin membrane. Here the two-scale homogenization asymptotics is performed simultaneously with a dimension reduction procedure -- a non-standard singular perturbation problem.  We numerically illustrate in Section \ref{numerics}  the behavior of the upscaled membrane in the finite thickness regime. Finally, in Section \ref{Discussion} we present our conclusions.

\section{The model}
\label{s:modello}
\noindent
Let $\ell,h>0$ and 
consider the two--dimensional strip 
$[-\ell/2,\ell/2]\times[0,h]$, say that 
$\ell$ and $h$ are, respectively, its horizontal and vertical 
side lengths. 
Partition the strip into the blocks
$\omega_\textrm{l}=[-\ell/2,-w/2]\times[0,h]$, 
$\omega_\textrm{m}=[-w/2,w/2]\times[0,h]$, 
$\omega_\textrm{r}=[w/2,\ell/2]\times[0,h]$, and 
call $\omega_\textrm{m}$ the \textit{membrane}. 
Let $0<\eta\le h$ and $\varepsilon=2\eta/\ell$. We 
partition the membrane into 
rectangular \textit{cells} 
$\omega_\textrm{c}^i=(-w/2,w/2)\times((i-1)\eta,i\eta)\cap(0,h)$ 
with $i$ running from one to the smallest integer larger than or 
equal to $h/\eta$. 
In each cell consider an impenetrable disk, 
called \textit{obstacle}, 
with center in the center of the cell and diameter $O(\varepsilon)$ 
in the limit $\varepsilon\to0$.
Denote by $\omega_\textrm{o}$ the union of all 
the obstacles. 

We denote by 
$\gamma_\textrm{v}$
and
$\gamma_\textrm{h}$,
the vertical and horizontal boundaries of the strip, 
by $\gamma_\rr{o}$ the boundary of the 
obstacle region $\omega_\textrm{o}$
and 
by $\gamma_\textrm{i}$ the boundary of the region $\omega_\textrm{i}$
for $\textrm{i}=\textrm{l,m,r}$. 
The external normal direction 
to a closed
curve is denoted here by $n$.

\begin{figure}[h]
\begin{picture}(400,80)(-140,40)
\thicklines
\multiput(0,0)(0,100){2}{\line(1,0){200}}
\multiput(0,0)(200,0){2}{\line(0,1){100}}
\thinlines
\multiput(80,0)(40,0){2}{\line(0,1){100}}
\multiput(80,0)(0,16){7}{\line(1,0){40}}

\put(125,24){\vector(0,1){8}}
\put(125,24){\vector(0,-1){8}}
\put(127,22){${\scriptstyle \eta}$}

\put(205,50){\vector(0,1){50}}
\put(205,50){\vector(0,-1){50}}
\put(207,48){${\scriptstyle h}$}

\put(-10,-10){${\scriptstyle -\ell/2}$}
\put(70,-10){${\scriptstyle -w/2}$}
\put(110,-10){${\scriptstyle +w/2}$}
\put(190,-10){${\scriptstyle +\ell/2}$}

\put(15,80){${\scriptstyle \gamma_\textrm{v}\cap\gamma_\textrm{l}}$}
\put(15,80){\vector(-1,-1){10}}
\put(15,20){${\scriptstyle \gamma_\textrm{h}}$}
\put(20,15){\vector(1,-1){10}}
\put(162,7){${\scriptstyle \gamma_\textrm{v}\cap\gamma_\textrm{r}}$}
\put(185,10){\vector(1,1){10}}
\put(177,80){${\scriptstyle \gamma_\textrm{h}}$}
\put(175,85){\vector(-1,1){10}}

\put(40,40){${\scriptstyle \omega_\textrm{l}}$}
\put(160,40){${\scriptstyle \omega_\textrm{r}}$}

\put(55,60){${\scriptstyle \omega_\textrm{m}}$}
\put(65,57){\vector(1,-1){20}}
\put(65,67){\vector(1,1){20}}

\put(136,61){${\scriptstyle \omega_\textrm{o}}$}
\put(135,64.5){\vector(-3,2){32}}
\put(135,62.5){\vector(-3,-2){32}}
\put(135,60){\vector(-2,-3){33}}

\linethickness{1.5mm}
\multiput(95,7.5)(0,16){6}{\line(1,0){10}}
\end{picture}
\vskip 2. cm 
\caption{Schematic representation of the model geometry.}
\label{f:mod000} 
\end{figure}
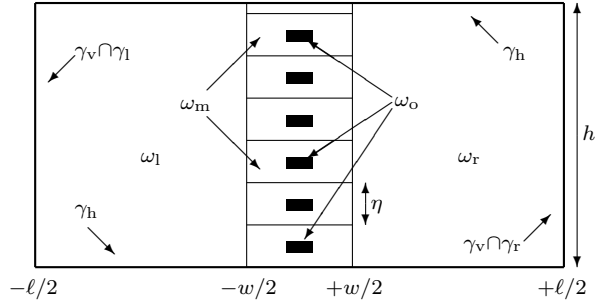


\noindent We let
$\omega
 =(\omega_\textrm{l}\cup\omega_\textrm{r}\cup
  \omega_\textrm{m})\setminus\omega_\textrm{o}$ and $f:\omega\to\bb{R}$ be a 
real function. Fixing the parameters $d_1,d_2>0$,
we consider the differential problem 
\begin{equation}
\label{mod000}
\frac{\partial u}{\partial t}
-d_1 \frac{\partial^2u}{\partial x_1^2}
-d_2 \frac{\partial^2u}{\partial x_2^2}
= -\frac{\partial}{\partial x_1}bu(1-u)+f(x)\quad \mbox{in}\,\,\omega,
\end{equation}
endowed 
with the homogeneous Neumann boundary conditions 
\begin{equation}
\label{mod005}
\Big(
d_1\frac{\partial u}{\partial x_1}-bu(1-u),
d_2\frac{\partial u}{\partial x_2}
\Big)
\cdot n=0
\;\;\;\textrm{ on } \gamma_\textrm{h}\cup\gamma_\textrm{o},
\end{equation}
as well as with the Dirichlet conditions
\begin{equation}
\label{mod010}
u(x,t)=u_\textrm{l}
\;\;\textrm{ on }
\gamma_\textrm{v}\cap\gamma_\textrm{l}
\;\;\textrm{ and }\;\;
u(x,t)=u_\textrm{r}
\;\;\textrm{ on }
\gamma_\textrm{v}\cap\gamma_\textrm{r}
\end{equation}
for any $t\ge0$, where $u_\textrm{l},u_\textrm{r}\in\bb{R}$. As initial condition we take 
\begin{equation}
\label{mod020}
u(x,0)=v(x)
\;\;\textrm{ on } \omega.
\end{equation}
\subsection{The non--dimensional model}
\label{s:adim}
\noindent
It is useful to introduce dimensionless variables 
\begin{equation}
\label{adim000}
X=(X_1,X_2)=
\Big(\frac{2x_1}{\ell},\frac{2x_2}{\ell}\Big)
\;\;\textrm{ and }\;\;
T=\frac{t}{\tau},\;
\end{equation}
where $\tau$ is a fixed positive real.

Using \eqref{adim000}, the original strip is mapped to  
$[-1,1]\times[0,2h/\ell]$, 
which is partitioned into 
$\Omega_\textrm{l}=[-1,-w/\ell]\times[0,2h/\ell]$,
$\Omega_\textrm{m}=[-w/\ell,w/\ell]\times[0,2h/\ell]$,
and
$\Omega_\textrm{r}=[w/\ell,1]\times[0,2h/\ell]$.
The cells are mapped to 
$\Omega_\textrm{c}^i=(-w/\ell,w/\ell)\times((i-1)\varepsilon,
                       i\varepsilon)\cap(0,2h/\ell)$, 
where we recall that $\varepsilon=2\eta/\ell$.
In the new variables, we denote by $\Omega_\textrm{o}$ the region 
occupied by the obstacle and by 
$\Gamma_\textrm{v}$,
$\Gamma_\textrm{h}$,
$\Gamma_\textrm{l}$,
$\Gamma_\textrm{m}$,
$\Gamma_\textrm{r}$,
and
$\Gamma_\textrm{o}$
the boundaries introduced above.


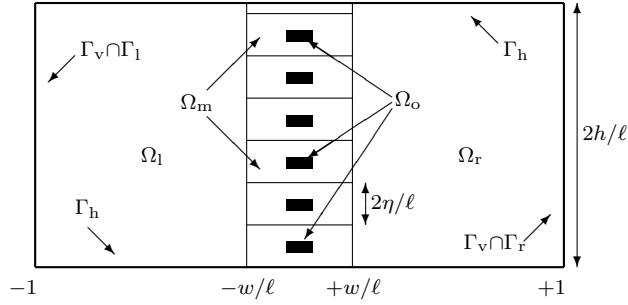
\begin{figure}[h]
\begin{picture}(400,80)(-140,40)
\thicklines
\multiput(0,0)(0,100){2}{\line(1,0){200}}
\multiput(0,0)(200,0){2}{\line(0,1){100}}
\thinlines
\multiput(80,0)(40,0){2}{\line(0,1){100}}
\multiput(80,0)(0,16){7}{\line(1,0){40}}

\put(125,24){\vector(0,1){8}}
\put(125,24){\vector(0,-1){8}}
\put(127,22){${\scriptstyle 2\eta/\ell}$}

\put(205,50){\vector(0,1){50}}
\put(205,50){\vector(0,-1){50}}
\put(207,48){${\scriptstyle 2h/\ell}$}

\put(-10,-10){${\scriptstyle -1}$}
\put(70,-10){${\scriptstyle -w/\ell}$}
\put(110,-10){${\scriptstyle +w/\ell}$}
\put(190,-10){${\scriptstyle +1}$}

\put(17,80){${\scriptstyle \Gamma_\textrm{v}\cap\Gamma_\textrm{l}}$}
\put(15,80){\vector(-1,-1){10}}
\put(15,20){${\scriptstyle \Gamma_\textrm{h}}$}
\put(20,15){\vector(1,-1){10}}
\put(162,7){${\scriptstyle \Gamma_\textrm{v}\cap\Gamma_\textrm{r}}$}
\put(185,10){\vector(1,1){10}}
\put(177,80){${\scriptstyle \Gamma_\textrm{h}}$}
\put(175,85){\vector(-1,1){10}}

\put(40,40){${\scriptstyle \Omega_\textrm{l}}$}
\put(160,40){${\scriptstyle \Omega_\textrm{r}}$}

\put(55,60){${\scriptstyle \Omega_\textrm{m}}$}
\put(65,57){\vector(1,-1){20}}
\put(65,67){\vector(1,1){20}}

\put(136,61){${\scriptstyle \Omega_\textrm{o}}$}
\put(135,64.5){\vector(-3,2){32}}
\put(135,62.5){\vector(-3,-2){32}}
\put(135,60){\vector(-2,-3){33}}

\linethickness{1.5mm}
\multiput(95,7.5)(0,16){6}{\line(1,0){10}}
\end{picture}
\vskip 2. cm 
\caption{Schematic representation of the dimensionless model geometry.}
\label{f:mod010} 
\end{figure}

It is convenient to set 
\begin{equation}
\label{adim005}
U(X,T)=u\Big(\frac{\ell X}{2}, \tau T\Big),\;
V(X)=v\Big(\frac{\ell X}{2}\Big),\;
F(X)=\tau f\Big(\frac{\ell X_1}{2},\frac{\ell X_2}{2}\Big)
\end{equation}
and rewrite 
the model \eqref{mod000} as follows
\begin{equation}
\label{adim010}
\frac{\partial U}{\partial T}
+\nabla\cdot J=F
\end{equation}
in 
$\Omega
 =(\Omega_\textrm{l}\cup\Omega_\textrm{r}\cup
  \Omega_\textrm{m})\setminus\Omega_\textrm{o}$,
where we introduced the flux 
\begin{equation}
\label{adim015}
J=-D(\nabla U+G(U))
\;,
\end{equation}
with
the derivatives in $\nabla$ taken  with respect to the dimensionless variables $X_1, X_2$,
and let 
\begin{equation}
\label{adim040}
D=\left(\begin{array}{cc}D_1&0\\0&D_2\end{array}\right),
\;
D_1=\frac{4\tau d_1}{\ell^2},\;
D_2=\frac{4\tau d_2}{\ell^2},
\;\textrm{ and }\;
G(U)=\left(\begin{array}{c}g(U)\\ 0\end{array}\right), \;
\;
\end{equation}
with $g(U)=\ell p(U)/(2d_1)$, where $p(U)=-bU(1-U)$ -- a  choice that makes (\ref{adim015})  to correspond precisely to the setting discussed in \cite{CKMSpre2016}. 


The derivations done in this paper cover the  more general case:
 \begin{equation}
\label{adim040-gen}
D=\left(\begin{array}{cc}D_{11}& D_{12} \\ D_{21}&D_{22}\end{array}\right),
\;
\;G(U)=\left(\begin{array}{c}g(U)\\ 0\end{array}\right), \;
\; \textrm{ and }\;
p(U) \textrm{ an arbitrary polynomial}.
\end{equation}
To fix ideas, we take $p(U)=\sum_{n=1}^ka_nU^n$, where $a_n\in\mathbb{R}$.
If not mentioned otherwise, in the rest of the paper $D$ is a full matrix as indicated  in  (\ref{adim040-gen}). 

For any
$T\ge0$, problem \eqref{adim010} is endowed with the Dirichlet boundary 
conditions
\begin{equation}
\label{adim020}
U(X,T)
=u_\textrm{l}
\;\;\textrm{ on }
\Gamma_\textrm{v}\cap\Gamma_\textrm{l}
\;\;\textrm{ and }\;\;
U(X,T)
=u_\textrm{r}
\;\;\textrm{ on }
\Gamma_\textrm{v}\cap\Gamma_\textrm{r},
\end{equation}
the 
Neumann boundary conditions
\begin{equation}\label{flux}
J\cdot n=0\quad \mbox{on}\,\,\Gamma_\textrm{h}\cup \Gamma_\textrm{o},
\end{equation}
and the initial condition 
\begin{equation}
\label{adim030}
U(X,0)=V(X)
\;\;\textrm{ in } \Omega.
\end{equation}

\section{Derivation of the finite-thickness upscaled membrane  model}
\label{s:duescale}
\par\noindent

In this section, we use a two-scale homogenization approach to average the membrane's internal structure and then to derive the corresponding upscaled equation for the mass transport 
as well as the effective transport coefficient. If the diffusion matrix is diagonal, then we point out explicitly the structure of the corresponding 
tortuosity tensor.

\subsection{Two-scale expansions}

We look for effective equations in the limit in which the height of the 
cells tends to zero and its number is increased so that the total 
height of the cells equals that of the whole strip. 
Due to the periodic micro--structure of the membrane $\Omega_\textrm{m}$, 
with vertical spatial period $\varepsilon=2\eta/\ell$, it is 
reasonable to attack the problem expanding the unknown function 
$U$ in the membrane region as  
\begin{equation}
\label{due010}
U(X,T)
=
\sum_{n=0}^\infty \varepsilon^n U^\textrm{m}_n(X,Y_2,T)
\;\;\;\;\;\;\;\textrm{in } \Omega_\textrm{m},
\end{equation}
where $Y_2=X_2/\varepsilon$ and the functions $U^\textrm{m}_n$ 
are $Y_2-$periodic functions.

By abusing slightly the notation, we understand in \eqref{adim010}
\begin{displaymath}
\nabla=\nabla_X+\frac{1}{\varepsilon}\nabla_{Y_2}
\;\;\textrm{with }
\nabla_X=\left(\begin{array}{c}
             \frac{\partial}{\partial X_1}\\
             \frac{\partial}{\partial X_2}\end{array}\right)
\;\;\textrm{and }
\nabla_{Y_2}=\left(\begin{array}{c}
             0\\
             \frac{\partial}{\partial Y_2}\end{array}\right)
\;.
\end{displaymath}
We now compute the various terms appearing in \eqref{adim010} in 
the different regions of $\Omega$.
We have
\begin{equation}
\label{due020}
\frac{\partial U}{\partial T}
=
\sum_{n=0}^\infty
\varepsilon^n\frac{\partial U^\textrm{m}_n}{\partial T}
\;\;\textrm{ and }\;\;
\frac{\partial U}{\partial X_1}
=
\sum_{n=0}^\infty
\varepsilon^n\frac{\partial U^\textrm{m}_n}{\partial X_1}
\;\;\;\;\;\;\;\textrm{ in } \Omega_\textrm{m}.
\end{equation}

For handling the terms involving the gradient $\nabla$, we have to distinguish 
the regions $\Omega_\textrm{l}$, $\Omega_\textrm{m}$,  
and $\Omega_\textrm{r}$. 
In $\Omega_\textrm{l}$ and $\Omega_\textrm{r}$ 
we simply have  $\nabla U(X,T)=\nabla U_0^l(X,T)$ in $\Omega_l$ and $\nabla U(X,T)=\nabla U_0^r(X,T)$ in $\Omega_r$. Instead of $\nabla U_0^l$ and $\nabla U_0^r$, we will use $\nabla U^l$ and $\nabla U^r$, respectively.

In $\Omega_\textrm{m}$, the computation of the gradient reads
\begin{equation}
\label{nablaUm}
\begin{split}
\nabla U
=
&
\nabla\sum_{n=0}^\infty\varepsilon^nU^\textrm{m}_n
=
\sum_{n=0}^\infty\varepsilon^n\nabla_XU^\textrm{m}_n
+
\sum_{n=0}^\infty\varepsilon^n\frac{1}{\varepsilon}\nabla_{Y_2}U^\textrm{m}_n
\\
=
&
\frac{1}{\varepsilon}\nabla_{Y_2}U^\textrm{m}_0
+
\sum_{n=0}^\infty\varepsilon^n(\nabla_XU^\textrm{m}_n
+\nabla_{Y_2}U^\textrm{m}_{n+1}).
\\
\end{split}
\end{equation}
Hence, it yields 
\begin{equation}
\label{due050}
\begin{array}{rcl}
\nabla\cdot D\nabla U
&\!\!=&\!\!
{\displaystyle
\frac{1}{\varepsilon^2} \nabla_{Y_2}\cdot D\nabla_{Y_2} U^\textrm{m}_0
+
\frac{1}{\varepsilon} \nabla_X\cdot D\nabla_{Y_2} U^\textrm{m}_0
}
\\
&\!\!&\!\!
{\displaystyle
+\sum_{n=0}^\infty
  \varepsilon^n
  \Big[
       \nabla_X\cdot D\nabla_X U^\textrm{m}_n
       +
       \nabla_X\cdot D\nabla_{Y_2} U^\textrm{m}_{n+1}
       +
       \frac{1}{\varepsilon} \nabla_{Y_2}\cdot D\nabla_X U^\textrm{m}_n
}
\\
&\!\!&\!\!
{\displaystyle
 \phantom{+\sum_{n=0}^\infty \varepsilon^n\Big[}
       +
       \frac{1}{\varepsilon} \nabla_{Y_2}\cdot D\nabla_{Y_2} U^\textrm{m}_{n+1}
  \Big]
}
\\
&\!\!=&\!\!
{\displaystyle
\frac{1}{\varepsilon^2} \nabla_{Y_2}\cdot D\nabla_{Y_2} U^\textrm{m}_0
+
\frac{1}{\varepsilon} 
 [
      \nabla_X\cdot D\nabla_{Y_2} U^\textrm{m}_0
      +
      \nabla_{Y_2}\cdot D\nabla_X U^\textrm{m}_0
}
\\
&\!\!&\!\!
{\displaystyle
 \phantom{\frac{1}{\varepsilon^2} \nabla_{Y_2}\cdot 
          D\nabla_{Y_2} U^\textrm{m}_0+ \frac{1}{\varepsilon} [
         }
      +
      \nabla_{Y_2}\cdot D\nabla_{Y_2} U^\textrm{m}_1
 ]
}
\\
&\!\!&\!\!
{\displaystyle
+\sum_{n=0}^\infty
  \varepsilon^n
  [
       \nabla_X\cdot D\nabla_X U^\textrm{m}_n
       +
       \nabla_X\cdot D\nabla_{Y_2} U^\textrm{m}_{n+1}
       +
       \nabla_{Y_2}\cdot D\nabla_X U^\textrm{m}_{n+1}
}
\\
&\!\!&\!\!
{\displaystyle
    \phantom{+\sum_{n=0}^\infty\varepsilon^n[}
       +
       \nabla_{Y_2}\cdot D\nabla_{Y_2} U^\textrm{m}_{n+2}
  ].
}
\end{array}
\end{equation}
Moreover, 
 we have  
\begin{equation}
\label{DG}
DG(U)
=\left(\begin{array}{c} D_1g(U)\\ 0\\ \end{array}\right)
=
DG(U^\textrm{m}_0)
+
\varepsilon
D
\left(\begin{array}{c} 
U^\textrm{m}_1\sum_{n=1}^k n b_n (U^\textrm{m}_0)^{n-1} \\ 
0\\ \end{array}\right)
+o(\varepsilon).
\end{equation}
It is worth noting already at this stage that if the matrix $D$ is diagonal, then (\ref{DG}) reduces to
\begin{equation}
\label{nablaDG}
\nabla\cdot DG(U)
=
\nabla_X\cdot DG(U^\textrm{m}_0)
+o(1)
\;\;.
\end{equation}

We consider now the equation inside the membrane region $\Omega_\textrm{m}$
at the lowest order $\varepsilon^{-2}$ and we find 
\begin{equation}
\label{eqm000}
\nabla_{Y_2}\cdot D\nabla_{Y_2} U^\textrm{m}_0
=
0.
\end{equation}
By expanding $J$ and by collecting the lowest $\varepsilon$ order, 
we get the Neumann boundary condition 
\begin{equation}
\label{eqm005}
(-D\nabla_{Y_2} U^\textrm{m}_0)\cdot n=0\quad \mbox{on}\,\, 
(\Gamma_\textrm{o}\cup \Gamma_\textrm{h})\cap \Omega_\textrm{m}
\end{equation}
and the following transmission boundary conditions: 
$$U^\textrm{m}_0(X,T)=U^\textrm{l}_0(X,T) \quad \mbox{on}\,\,
\Gamma_\textrm{l}\cap \Gamma_\textrm{m}
\;\;\textrm{ and }\;\;
U^\textrm{m}_0(X,T)=U^\textrm{r}_0(X,T)\quad \mbox{on}\,\,
\Gamma_\textrm{r}\cap \Gamma_\textrm{m}$$
as well as
 \begin{equation}
  -D(\nabla U^\textrm{l}+G(U^\textrm{l}))\cdot n =-D(\nabla U^\textrm{m}_0+G(U^\textrm{m}_0)) \mbox{ at } \Gamma_\textrm{l}\cap \Gamma_\textrm{m},
 \end{equation}
  \begin{equation}
  -D(\nabla U^\textrm{r}+G(U^\textrm{r}))\cdot n=-D(\nabla U^\textrm{m}_0+G(U^\textrm{m}_0)) \mbox{ at } \Gamma_\textrm{r}\cap \Gamma_\textrm{m},
 \end{equation}
for any $T\ge0$.

We recall that $U^\textrm{m}_0$ is $Y_2-$periodic.
 Based on \eqref{eqm000} and \eqref{eqm005}, we claim that 
 $U^\textrm{m}_0$ is independent of $Y_2$, i.e. 
$U^\textrm{m}_0=U^\textrm{m}_0(X,T)$.

At the order $\varepsilon^{-1}$,
using that $U^\textrm{m}_0$ does not depend on $Y_2$, we get the equation
\begin{equation}
\label{eqm010}
 \nabla _{Y_2}\cdot D\nabla_{Y_2}U^\textrm{m}_1=-\nabla_{Y_2}\cdot D\nabla_XU^\textrm{m}_0
\end{equation}
with Neumann boundary condition \eqref{flux} at order $\varepsilon^0$ in \eqref{nablaUm} and \eqref{DG}  
\begin{equation}
\label{boundaryMembrana-1}
-D\nabla_{Y_2}U^\textrm{m}_1\cdot n=D\nabla_X U^\textrm{m}_0\cdot n
+DG(U^\textrm{m}_0)\cdot n \mbox{ on } \Gamma_h\cup \Gamma_o.
\end{equation}

Recall that $U^\textrm{m}_1$ is $Y_2-$periodic.

\subsection{$D$ diagonal matrix}

If $D$ is a diagonal matrix, then the structure of  \eqref{eqm010} allows us to assume that
\begin{equation}\label{W}
U^\textrm{m}_1=W(Y_2)\cdot \left[\nabla_XU^\textrm{m}_0
+G(U^\textrm{m}_0)\right],
\end{equation}
where $W(Y_2)$ is a vector with $Y_2-$periodic components. We will refer to  $W(Y_2)$ as {\em cell function}. 
Substituting now the expression \eqref{W} in \eqref{eqm010}, we get
$$
\nabla_{Y_2}\cdot D\nabla_{Y_2}W(Y_2)\cdot\left[\nabla_X U^\textrm{m}_0
+G(U^\textrm{m}_0)\right]
=-\nabla_{Y_2}\cdot D\cdot\left[\nabla_X U^\textrm{m}_0
+G(U^\textrm{m}_0)\right],
$$
while substituting the same expression now in \eqref{boundaryMembrana-1} leads to
 $$ -D\nabla_{Y_2}W(Y_2)\cdot\left[\nabla_X U^\textrm{m}_0
+G(U^\textrm{m}_0)\right]\cdot n = D\left[\nabla_X U^\textrm{m}_0
+G(U^\textrm{m}_0)\right]\cdot n.$$

Now, we can introduce the following {\em cell problems}: find the $Y_2$-periodic cell function $W=(w_1,w_2)^T$ satisfying the   following elliptic partial differential equations:
 \begin{eqnarray}
\nabla_{Y_2}\cdot (D\nabla_{Y_2}w_j(Y_2))=-\nabla_{Y_2}\cdot D e_j, \label{cell1}\\
\nabla_{Y_2}w_j\cdot n=0 \mbox{ on } \Gamma_h\cup \Gamma_o,\label{cell2}
 \end{eqnarray}
 for $j=1,2.$
 In \eqref{cell1}, we use the coordinate vectors $e_1=(1 \ \ 0)^T$ and $e_2=(0 \ \ 1)^T$. We point out that \eqref{cell1} can be written explicitly as $\frac{\partial}{\partial Y_2}\left( D_{22}\frac{\partial w_1}{\partial Y_2}\right)=0$ and $\frac{\partial}{\partial Y_2}\left[ D_{22}\left(1+\frac{\partial w_1}{\partial Y_2}\right)\right]=0$, which in the absence of the internal heterogeneity can be solved analytically; see Proposition 3.3, p. 13 in \cite{Hornung}. 


For $U^\textrm{m}_2$, taking into account \eqref{due020}, 
\eqref{due050}, and \eqref{nablaDG}, at the order $\varepsilon^0$, 
we have the following equation
\begin{equation}
\label{complete}
\begin{split}
\frac{\partial U^\textrm{m}_0}{\partial T}
-[
&
\nabla_X\cdot D\nabla_XU^\textrm{m}_0
+\nabla_X\cdot D\nabla_{Y_2}U^\textrm{m}_1
\\
&
+\nabla_{Y_2}\cdot D\nabla_{X}U^\textrm{m}_1
+\nabla_{Y_2}\cdot D\nabla_{Y_2}U^\textrm{m}_2
+\nabla_X\cdot DG(U^\textrm{m}_0)]=F
\end{split}
\end{equation}
satisfying as boundary condition \eqref{flux} across $(\Gamma_\textrm{o}\cup \Gamma_\textrm{h})\cap \Omega_\textrm{m}$
\begin{equation}\label{bcU2}
-D\left[\nabla_XU^\textrm{m}_1+\nabla_{Y_2}U^\textrm{m}_2
+\left(\begin{array}{c}U_1^\textrm{m}
\sum_{n=0}^kb_nn(U^\textrm{m}_0)^{n-1}\\0\end{array}\right)  \right]\cdot n=0,
\end{equation}
obtained by using the order $\varepsilon$ of the 
expansions \eqref{nablaUm} and \eqref{DG}.


Integrating \eqref{complete} with respect to $Y_2$ over a cell, say on
the set $Z=[0,2\eta/\ell]$, using 
the divergence theorem with respect to the variable $Y_2$ 
and \eqref{W}, we have
\begin{displaymath}
\begin{split}
\int_Z \frac{\partial U^\textrm{m}_0}{\partial T}\textrm{d}Y_2
&
-\nabla_X\cdot\!\int_Z D\nabla_XU^\textrm{m}_0\textrm{d}Y_2
-\nabla_X\cdot\!\int_Z D\nabla_{Y_2}\left[W(Y_2)\cdot\left(\nabla_X U_0^\textrm{m}
+G(U^\textrm{m}_0)\right)\right]
\textrm{d}Y_2
\\
&
-\nabla_X\cdot \!\int_Z DG(U^\textrm{m}_0)\textrm{d}Y_2
-\!\int_Z\nabla_{Y_2}\cdot D\nabla_XU_1^\textrm{m}\textrm{d}Y_2
=\!\int_ZF\textrm{d}Y_2
+\!\int_{\partial Z}D\nabla_{Y_2}U_2^\textrm{m}\cdot n\textrm{d}\sigma.
\end{split}
\end{displaymath}
Notice that the last term in the above 
equation is noting but the differences between the values 
of the function $D\nabla_{Y_2}U_2^\textrm{m}\cdot n$ evaluated 
at the extremes $2\eta/\ell$ and $0$ of the integration interval. 
In that term $n$ is the external normal to the horizontal parts 
of the boundary of the elementary cell, in particular 
it is a vertical unit vector.
Hence, by using \eqref{bcU2}, we obtain  
\begin{displaymath}
\begin{split}
\int_Z \frac{\partial U^\textrm{m}_0}{\partial T}\textrm{d}Y_2
&
-\nabla_X\cdot\!\int_Z D\nabla_XU^\textrm{m}_0\textrm{d}Y_2
-\nabla_X\cdot\!\int_Z D\nabla_{Y_2}\left[W(Y_2)\cdot\left(\nabla_X U_0^\textrm{m}
+G(U^\textrm{m}_0)\right)\right]
\textrm{d}Y_2
\\
&
-\nabla_X\cdot \!\int_Z DG(U^\textrm{m}_0)\textrm{d}Y_2
-\!\int_Z\nabla_{Y_2}\cdot D\nabla_XU_1^\textrm{m}\textrm{d}Y_2
=\!\int_ZF\textrm{d}Y_2
-\!\int_{\partial Z}D\nabla_XU_1^\textrm{m}\cdot n\textrm{d}\sigma.
\end{split}
\end{displaymath}
By the divergence theorem, the last term of the left--hand side cancels 
the last term of the right--hand side. Thus, we get
\begin{displaymath}
\begin{split}
\int_Z \frac{\partial U^\textrm{m}_0}{\partial T}\textrm{d}Y_2
&
-\nabla_X\cdot\int_Z D[\nabla_XU^\textrm{m}_0+G(U^\textrm{m}_0)]\textrm{d}Y_2
-\nabla_X\cdot\int_Z D\nabla_{Y_2}\left[W(Y_2)\cdot\left(\nabla_X U_0^\textrm{m}
+G(U^\textrm{m}_0)\right)\right]
\textrm{d}Y_2
\\
&
=\int_ZF\textrm{d}Y_2
\;\;.
\end{split}
\end{displaymath}
Recalling that $U^\textrm{m}_0$ does not depend on $Y_2$, we finally get 
\begin{equation}\label{final}
\frac{\partial U^\textrm{m}_0}{\partial T}
-\nabla_X\cdot 
\left[\frac{1}{|Z|}\int_ZD\left(\mathbb{I}+\left(\begin{array}{cc} 0 & 0 \\ \frac{\partial w_1}{\partial Y_2}& \frac{\partial w_2}{\partial Y_2}\end{array} \right )\right)\textrm{d}Y_2\right]
\left(\nabla_X U^\textrm{m}_0
+G(U^\textrm{m}_0)\right)
=\frac{1}{|Z|}\int_Z F\textrm{d}Y_2.
\end{equation}
We refer to the coefficient 
\begin{equation}
\mathbb{D}:=\frac{1}{|Z|}\int_ZD\left(\mathbb{I}+\left(\begin{array}{cc} 0 & 0 \\ \frac{\partial w_1}{\partial Y_2}& \frac{\partial w_2}{\partial Y_2}\end{array}\right)\right)\textrm{d}Y_2
\end{equation}
 as effective transport coefficient. 

The upscaled equation \eqref{final} for the zero term of the 
expansion has the same structure as the 
original equation \eqref{adim010}. 
The source term $F$ on the right--hand side 
is replaced by its average 
over the cell on the $Y_2$.  
The diffusion matrix 
is replaced by 
its average over the cell on the $Y_2$ variable weighted by the function 
$$\mathbb{I}+\left(\begin{array}{cc} 0 & 0 \\ \frac{\partial w_1}{\partial Y_2}& \frac{\partial w_2}{\partial Y_2}\end{array}\right),$$ which is referred to as tortuosity tensor in the porous media literature; we refer the reader to the review paper \cite{Ekeoma} for a discussion done in terms of this tortuosity tensor of the role played by microscopic anisotropies in understanding macroscopically a smoldering combustion scenario.

Summarizing, the upscaled model equation reads:

Find $U_0^\textrm{m}(X,Y_1,T)$ satisfying
\begin{equation}\label{final01}
\frac{\partial U^\textrm{m}_0}{\partial T}
-\nabla_X\cdot 
\left[\frac{1}{|Z|}\int_ZD\left(\mathbb{I}+\left(\begin{array}{cc} 0 & 0 \\ \frac{\partial w_1}{\partial Y_2}& \frac{\partial w_2}{\partial Y_2}\end{array} \right )\right)\textrm{d}Y_2\right] 
\left(\nabla_X U^\textrm{m}_0
+G(U^\textrm{m}_0)\right)
=\frac{1}{|Z|}\int_Z F\textrm{d}Y_2.
\end{equation}
 \begin{equation}
 U_0^\textrm{m}=U^\textrm{l}, -D(\nabla U^\textrm{l}+G(U^\textrm{l}))\cdot n=-D(\nabla_XU^\textrm{m}_0+G(U^\textrm{m}_0))\cdot n \mbox{ at } \Gamma_\textrm{l}\cap \Gamma_\textrm{m},
 \end{equation}
  \begin{equation}
 U_0^\textrm{m}=U^\textrm{r}, -D(\nabla U^\textrm{r}+G(U^\textrm{r}))\cdot n=-D(\nabla_XU^\textrm{m}_0+G(U^\textrm{m}_0))\cdot n \mbox{ at } \Gamma_\textrm{r}\cap \Gamma_\textrm{m},
 \end{equation}
 together with the initial condition
 \begin{equation}\label{final0f}
 U_0^\textrm{m}(T=0)=V^\textrm{m}(X,Y_1).
 \end{equation}
Using the transmission conditions at $\Gamma_\textrm{l}$ and $\Gamma_\textrm{r}$, the information in $\Omega_\textrm{m}$ is now linked (in a well-posed way) with equation \eqref{adim010} posed in $\Omega_\textrm{l}$ and $\Omega_\textrm{r}$, respectively.

\subsection{$D$ full matrix}

If $D$ is a genuine full matrix, then $U_1^m$ cannot be expressed in a convenient closed form in terms of cell functions. In this case, the resulting upscaled system of equations reads: 

Find ($U_0^\textrm{m}(X,T), U_1^\textrm{m}(X,Y_1,Y_2,T)$) satisfying the following system of equations:

\begin{equation}
\frac{\partial U^\textrm{m}_0}{\partial T}
-\nabla_X\cdot \frac{1}{|Z|}\int_Z D[\nabla_XU^\textrm{m}_0+G(U^\textrm{m}_0)]\textrm{d}Y_2
-\nabla_X\cdot\frac{1}{|Z|}\int_Z D\nabla_{Y_2} U_1^m \textrm{d}Y_2
=\frac{1}{|Z|}\int_ZF\textrm{d}Y_2
\end{equation}
coupled with
\begin{equation}
 \nabla _{Y_2}\cdot D\nabla_{Y_2}U^\textrm{m}_1=-\nabla_{Y_2}\cdot D\nabla_XU^\textrm{m}_0,
\end{equation}
provided the following boundary conditions are given
\begin{equation}
-D\nabla_{Y_2}U^\textrm{m}_1\cdot n=D\nabla_X U^\textrm{m}_0\cdot n
+DG(U^\textrm{m}_0)\cdot n \mbox{ on } \Gamma_h\cup \Gamma_o,
\end{equation}
\begin{equation}
 U^\textrm{m}_1 \mbox{ is } Y_2-\mbox{periodic},
 \end{equation}
 \begin{equation}
 U_0^\textrm{m}=U^\textrm{l}, -D(\nabla U^\textrm{l}+G(U^\textrm{l}))\cdot n =-D(\nabla_XU^\textrm{m}_0+\nabla_{Y_2}U_1^\textrm{m}+G(U^\textrm{m}_0))\cdot n \mbox{ at } \Gamma_\textrm{l}\cap \Gamma_\textrm{m},
 \end{equation}
  \begin{equation}
 U_0^\textrm{m}=U^\textrm{r}, -D(\nabla U^\textrm{r}+G(U^\textrm{r}))\cdot n=-D(\nabla_XU^\textrm{m}_0+\nabla_{Y_2}U_1^m+G(U^\textrm{m}_0))\cdot n \mbox{ at } \Gamma_\textrm{r}\cap \Gamma_\textrm{m},
 \end{equation}
 together with the initial condition
 \begin{equation}
 U_0^\textrm{m}(T=0)=V^\textrm{m}(X,Y_1).
 \end{equation}
 As in the previous section, using the transmission conditions at $\Gamma_\textrm{l}$ and $\Gamma_\textrm{r}$, the information in $\Omega_\textrm{m}$ is now linked (in a well-posed way) with equation \eqref{adim010} posed in $\Omega_\textrm{l}$ and $\Omega_\textrm{r}$, respectively. 
 
\section{Derivation of the infinitely-thin upscaled membrane  model}
\label{s:infinitesima}
\par\noindent

We look for the effective model in the limit in which both the 
width and the height of the 
cells tends to zero and its number is increased so that the total 
height of the cells equals that of the whole strip.
In this limit the evolutive equation inside the membrane must be replaced 
by a matching condition between the solutions of the problems in 
the left and the right regions 
$\Omega_\textrm{l}$ and 
$\Omega_\textrm{r}$. In this case, the upscaling procedure needs to be combined with a singular perturbation ansatz; see \cite{Ene} for a remotely related case.

\subsection{Two-scale layer expansions}

We  consider the geometry introduced in 
Section~\ref{s:adim} and assume $w=2\eta$, so that 
the membrane is the region $[-2\eta/\ell,2\eta/\ell]\times[0,2h/\ell]$
(see Figure~\ref{f:mod010}). Recalling the relation $\varepsilon=2\eta/\ell$, 
in the homogenization limit $\varepsilon\to0$ the membrane shrinks to an 
infinitesimal wide separating surface.
The equations in 
$\Omega_\textrm{l}$ and $\Omega_\textrm{r}$ are as in Section~\ref{s:adim}, 
see equations \eqref{adim010}--\eqref{adim040}. More precisely, 
we have 
\begin{equation}
\label{inf000}
\frac{\partial U^\textrm{i}}{\partial T}
+\nabla\cdot J^\textrm{i}
=
F^\textrm{i}
\;\textrm{ in }\;\Omega_\textrm{i}
\;\textrm{ with }\;
J^\textrm{i}=-D^\textrm{i}(\nabla U^\textrm{i}+G(U^\textrm{i}))
\;\textrm{ for i}=\textrm{l,r},
\end{equation}
where $F^\textrm{l}:[-1,-\varepsilon]\to\mathbb{R}$,
$F^\textrm{r}:[\varepsilon,+1]\to\mathbb{R}$, $D^\textrm{i}$ a general real 
$2\times2$ matrix,
and 
\begin{equation}
\label{inf010}
G(U)=\left(\begin{array}{c}g(U)\\0\\\end{array}\right)
\end{equation}
with $g(U)=\sum_{n=1}^kb_nU^n$ where $b_n$ are real coefficients. 
In the membrane $\Omega_\textrm{m}\setminus\Omega_\textrm{o}$,  
we consider the equation 
\begin{equation}
\label{inf020}
\frac{1}{\varepsilon}
\frac{\partial U^\textrm{m}}{\partial T}
+\nabla\cdot J^\textrm{m}
=
\frac{1}{\varepsilon}
F^\textrm{m}\Big(\frac{X_1}{\varepsilon},X_2\Big)
\end{equation}
with $F^\textrm{m}:[-1,+1]\times[0,2h/\ell]\to\mathbb{R}$ and
the flux $J^\textrm{m}$ defined as 
\begin{equation}
\label{inf030}
J^\textrm{m}=-D^\textrm{m}\Big(\frac{X_1}{\varepsilon},X_2\Big)
(\varepsilon\nabla U^\textrm{m}+G(U^\textrm{m}))
\;,
\end{equation}
where $D^\textrm{m}$ is a $2\times2$ square matrix
$$D^\textrm{m}=\left(\begin{array}{cc}D^\textrm{m}_{11}&D^\textrm{m}_{12}\\D^\textrm{m}_{21}&D^\textrm{m}_{22}\end{array}\right).
\;
$$

These equations are endowed with 
the Dirichlet boundary 
conditions
\begin{equation}
\label{inf050}
U^\textrm{l}(X,T)
=u_\textrm{l}
\;\;\textrm{ on }
\Gamma_\textrm{v}\cap\Gamma_\textrm{l}
\;\;\textrm{ and }\;\;
U^\textrm{r}(X,T)
=u_\textrm{r}
\;\;\textrm{ on }
\Gamma_\textrm{v}\cap\Gamma_\textrm{r}
\end{equation}
for any
$T\ge0$, the initial condition 
\begin{equation}
\label{inf060}
U^\textrm{i}(X,0)=V^\textrm{i}(X)
\;\textrm{ in } \Omega_\textrm{i}\;
\;\textrm{ for i}=\textrm{l,r}\;
\;\textrm{ and }\;
U^\textrm{m}(X,0)=V^\textrm{m}(X)
\;\;\textrm{ in } \Omega_\textrm{m}\setminus\Omega_\textrm{o},
\end{equation}
the Neumann boundary conditions
\begin{equation}
\label{inf070}
J^\textrm{i}(X,T)\cdot n=0
\;\textrm{ on } \Gamma_\textrm{h}\cap\Omega_\textrm{i}\;
\;\textrm{ for i}=\textrm{l,r}\;
\;\textrm{ and }\;
J^\textrm{m}(X,T)\cdot n=0
\;\;\textrm{ on } 
(\Gamma_\textrm{h}\cap\Omega_\textrm{m})\cup\Gamma_\textrm{o}
\end{equation}
for any $T\ge0$,
the continuity (linear transmission) conditions 
\begin{equation}
\label{inf080}
U^\textrm{i}(X,T)=U^\textrm{m}(X,T) 
\;\textrm{ and }\;
J^\textrm{i}(X,T)\cdot n
=
J^\textrm{m}(X,T)\cdot n
\quad \mbox{on}\,\,
\Gamma_\textrm{i}\cap \Gamma_\textrm{m}
\;\textrm{ for i}=\textrm{l,r}\;
\end{equation}
for any $T\ge0$,
where in the last equation $n$ is the horizontal unit vector 
pointing to the left on 
$\Gamma_\textrm{l}$ 
and to the right on 
$\Gamma_\textrm{r}$.

Inside the membrane we use the same two--scale expansion as 
the one introduced in the Section~\ref{s:duescale}, namely we take 
\begin{equation}
\label{inf110}
U^\textrm{m}(X,T)
=
\sum_{n=0}^\infty \varepsilon^n U^\textrm{m}_n(X,y_2,T)
\;\;\;\textrm{in } \Omega_\textrm{m},
\end{equation}
where $y_2=X_2/\varepsilon$ and the functions $U^\textrm{i}_n$, 
with $\textrm{i}=\textrm{l,m,r}$, 
are $y_2-$periodic functions. Since the domain where the two-scale expansion is defined vanishes as $\varepsilon \to 0$, we refer to \eqref{inf110} as two--scale layer expansion. We claim that this expansion discovers formally precisely the limit point of the two-scale convergence for thin homogeneous layers (as presented cf. Definition 4.1 in \cite{membrane}). 

We define the new variables 
\begin{equation}
\label{inf200}
z_1=\frac{X_1}{\varepsilon}
\;\;\textrm{ and }\;\;
z_2=X_2,
\end{equation}
and, abusing the notation (recall, indeed, that small $u$ had a 
different meaning in Section~\ref{s:modello}), we set  
\begin{equation}
\label{inf210}
u^\textrm{m}(z,T)=U^\textrm{m}(\varepsilon z_1,z_2,T)
\end{equation}
for the original functions 
and 
\begin{equation}
\label{inf220}
u^\textrm{m}_n(z,y_2,T)=U^\textrm{m}_n(\varepsilon z_1,z_2,y_2,T)
\end{equation}
for the perturbative terms $n\ge0$. 

It is immediate to deduce the following derivation rules with respect 
to the new variables. We let 
\begin{equation}
\label{inf230}
\nabla_{z_1}
=
\left(\begin{array}{c}\frac{\partial}{\partial z_1}\\0\end{array}\right)
,\;\;
\nabla_{z_2}
=
\left(\begin{array}{c}0\\\frac{\partial}{\partial z_2}\end{array}\right)
,\;\;\textrm{ and }\;\;
\nabla_{y_2}
=
\left(\begin{array}{c}0\\\frac{\partial}{\partial y_2}\end{array}\right)
\end{equation}
and prove 
\begin{equation}
\label{inf250}
\nabla U^\textrm{m}_n
=
\frac{1}{\varepsilon}\nabla_{z_1}u^\textrm{m}_n
+\nabla_{z_2}u^\textrm{m}_n
+\frac{1}{\varepsilon}\nabla_{y_2}u^\textrm{m}_n
\;\;\textrm{ for } n\ge0.
\end{equation}
 
Firstly, we note that the first term $\varepsilon^0$ in the 
expansion of $J^\textrm{m}$ is 
\begin{equation}
\label{inf255}
J^\textrm{m}
=
-D^\textrm{m}\nabla_{z_1}u_0^\textrm{m}
-D^\textrm{m}\nabla_{y_2}u_0^\textrm{m}
-\left(\begin{array}{c}
D_{11}^\textrm{m}g(u^\textrm{m}_0)
\\
D_{21}^\textrm{m}g(u^\textrm{m}_0)
\end{array}\right)
+o(1)
\;\;.
\end{equation}
Hence, 
expanding the equation \eqref{inf020} in the region 
$\Omega_\textrm{m}\setminus\Omega_\textrm{o}$ and taking into account 
the order $\varepsilon^{-1}$
we get the following equation
\begin{equation}
\label{inf260}
\begin{split}
\frac{\partial u_0^\textrm{m}}{\partial T}
&
-\left[\nabla_{z_1}\cdot D^\textrm{m}\nabla _{z_1}u_0^\textrm{m}
 +\nabla_{y_2}\cdot D^\textrm{m}\nabla_{y_2}u_0^\textrm{m}
 +\nabla_{z_1}\cdot D^\textrm{m}\nabla_{y_2}u_0^\textrm{m}
 +\nabla_{y_2}\cdot D^\textrm{m}\nabla _{z_1}u_0^\textrm{m}\right]
\\
&
-\frac{\partial}{\partial z_1}(D^\textrm{m}_{11}g(u_0^\textrm{m}))
-\frac{\partial}{\partial y_2}(D^\textrm{m}_{21}g(u_0^\textrm{m}))
=F^{\textrm{m}}.
\end{split}
\end{equation}

We remark that in the limit $\varepsilon\to0$ the function $u^\textrm{m}_0$ 
will depend only on $T$, $z_2$, and $y_2$, that is to say 
the dependence on $z_1$ will be lost. One can see this effect in the last equation, if one rescales the variables back to $X_1=\varepsilon z_1$. Consequently, three terms will be proportional to $\varepsilon$.  Hence, the limit function $u_0^\textrm{m}$
will have to solve the equation 
\begin{displaymath}
\frac{\partial u_0^\textrm{m}}{\partial T}
-\nabla_{y_2}\cdot D^\textrm{m}\nabla_{y_2}u_0^\textrm{m}
-\frac{\partial}{\partial y_2}(D^\textrm{m}_{21}g(u_0^\textrm{m}))
=F^{\textrm{m}}, 
\end{displaymath}
which can be rewritten as 
\begin{equation}
\label{inf265}
\frac{\partial u_0^\textrm{m}}{\partial T}
-\nabla_{y_2}\cdot D^\textrm{m}[\nabla_{y_2}u_0^\textrm{m}
+ G(u_0^\textrm{m})]
=F^{\textrm{m}}
\end{equation}
for any $X_2$. 
The limit function $u^\textrm{m}_0$ is periodic in $y_2$ and 
has to satisfy the conditions
\begin{equation}
\label{inf267}
u^\textrm{m}_0(z_2,y_2,T)
=
U^\textrm{i}(0,z_2,T)
\textrm{ for i}=\textrm{l,r}\textrm{ and }
u^\textrm{m}_0(z_2,y_2,0)
=
V^\textrm{m}(X_1,z_2).
\end{equation}

In the limit $\varepsilon\to0$ the functions $U^\textrm{i}$, with
$\textrm{i}=\textrm{l,r}$ will solve the equations 
\eqref{inf000} with the conditions \eqref{inf050}, 
\eqref{inf060} (first equation), and \eqref{inf070} (first equation).
Moreover, the matching conditions \eqref{inf080} will provide as 
with a jump condition on the flux associated to the 
limit solutions $U^\textrm{i}$. Indeed, we first note that 
at order $\varepsilon^0$, using \eqref{inf255}, 
the matching condition \eqref{inf080} (second equation) can be 
written as 
\begin{equation}
\label{inf270}
-D^\textrm{l}\left(\nabla U^\textrm{l}+G(U^\textrm{l})\right)\cdot n
=
D^\textrm{m}_{11}\frac{\partial u^\textrm{m}_0}{\partial z_1}
+
D^\textrm{m}_{12}\frac{\partial u^\textrm{m}_0}{\partial y_2}
+
D^\textrm{m}_{11}g(u^\textrm{m}_0)
\end{equation}
and 
\begin{equation}
\label{inf275}
-D^\textrm{r}(\nabla U^\textrm{r}+G(U^\textrm{r}))\cdot n
=
-
D^\textrm{m}_{11}\frac{\partial u^\textrm{m}_0}{\partial z_1}
-
D^\textrm{m}_{12}\frac{\partial u^\textrm{m}_0}{\partial y_2}
-
D^\textrm{m}_{11}g(u^\textrm{m}_0).
\end{equation}
It is worth noting that equations (\ref{inf270}) and (\ref{inf275}) complete the system of upscaled equations; compare e.g. how Corollary 7.1  in \cite{membrane} proves a similar statement. These conditions emphasize that the macroscopic flux is obtained by averaging the corresponding microscopic flux.

\subsection{Summary of the upscaled equations}

The resulting upscaled problem corresponding to this asymptotic regime is:

Find the triplet $(U^\textrm{l}, u^\textrm{m}_0, U^\textrm{r})$ satisfying the following set of equations:
\begin{equation}\label{bi}
\frac{\partial U^\textrm{i}}{\partial T}
+\nabla\cdot [-D^\textrm{i}(\nabla U^\textrm{i}+G(U^\textrm{i}))]
=
F^\textrm{i}
\;\textrm{ in }\;\Omega_\textrm{i}, i=\textrm{l},\textrm{r},
\end{equation}

\begin{equation}
\frac{\partial u_0^\textrm{m}}{\partial T}
-\nabla_{y_2}\cdot D^\textrm{m}[\nabla_{y_2}u_0^\textrm{m}
+ G(u_0^\textrm{m})]
=F^{\textrm{m}},
\end{equation}

\begin{equation}
u^\textrm{m}_0 \mbox{ is periodic in } y_2
\end{equation}

\begin{equation}
u^\textrm{m}_0(z_2,y_2,T)
=
U^\textrm{i}(0,z_2,T)
\textrm{ for i}=\textrm{l,r}\textrm{ and }
u^\textrm{m}_0(z_2,y_2,0)
=
V^\textrm{m}(X_1,z_2)
\end{equation}
\begin{equation}
 -D^\textrm{l}(\nabla U^\textrm{l}+G(U^\textrm{l})\cdot n
=
D^\textrm{m}_{11}\frac{\partial u^\textrm{m}_0}{\partial z_1}
+
D^\textrm{m}_{12}\frac{\partial u^\textrm{m}_0}{\partial y_2}
+
D^\textrm{m}_{11}g(u^\textrm{m}_0),
\end{equation}
\begin{equation}
-D^\textrm{r}(\nabla U^\textrm{r}+G(U^\textrm{r}))\cdot n
=
-
D^\textrm{m}_{11}\frac{\partial u^\textrm{m}_0}{\partial z_1}
-
D^\textrm{m}_{12}\frac{\partial u^\textrm{m}_0}{\partial y_2}
-
D^\textrm{m}_{11}g(u^\textrm{m}_0),
\end{equation}

\begin{equation}
U^\textrm{l}(X,T)
=u_\textrm{l}
\;\;\textrm{ on }
\Gamma_\textrm{v}\cap\Gamma_\textrm{l}
\;\;\textrm{ and }\;\;
U^\textrm{r}(X,T)
=u_\textrm{r}
\;\;\textrm{ on }
\Gamma_\textrm{v}\cap\Gamma_\textrm{r},
\end{equation}

\begin{equation}
J^\textrm{i}(X,T)\cdot n=0
\;\textrm{ on } \Gamma_\textrm{h}\cap\Omega_\textrm{i}\;
\;\textrm{ for i}=\textrm{l,r},
\end{equation}
\begin{equation}\label{bf}
U^\textrm{i}(X,0)=V^\textrm{i}(X)
\;\textrm{ in } \Omega_\textrm{i}\;
\;\textrm{ for i}=\textrm{l,r}.
\end{equation}

\subsection{Further remarks}

In what follows, we deduce alternative transmission relations across the membrane, recovering expected structures as if one would have applied  two-scale layer convergence arguments as indicated in \cite{membrane}.
 
Integrating the equation \eqref{inf260} with respect to $z_1$ we get
\begin{displaymath}
\begin{split}
\int_{-1}^{1}
&
\frac{\partial u_0^\textrm{m}}{\partial T}\textrm{d}z_1
-\left[D^\textrm{m}_{11}\frac{\partial u^\textrm{m}_0}{\partial z_1}
 \right]_{z_1=-1}^{z_1=+1}
-\nabla_{y_2}\cdot\int_{-1}^{1}\!\!\!D^\textrm{m}\nabla_{y_2}u_0^\textrm{m}
 \textrm{d}z_1
-\left[D^\textrm{m}_{12}\frac{\partial u^\textrm{m}_0}{\partial y_2}
 \right]_{z_1=-1}^{z_1=+1}
\\
&
-\nabla_{y_2}\cdot\int_{-1}^{1}\!\!\!D^\textrm{m}\nabla _{z_1}u_0^\textrm{m}
 \textrm{d}z_1
-\left[D_{11}^\textrm{m}g(u_0^\textrm{m})\right]_{z_1=-1}^{z_1=1}
-\int_{-1}^{1}\frac{\partial}{\partial y_2}(D^\textrm{m}_{21}g(u_0^\textrm{m}))
 \textrm{d}z_1
=\int_{-1}^{1}F^\textrm{m}\textrm{d}z_1
.
\end{split}
\end{displaymath}
By \eqref{inf270} and \eqref{inf275} we get
\begin{displaymath}
\begin{split}
\int_{-1}^{1}
&
\frac{\partial u_0^\textrm{m}}{\partial T}\textrm{d}z_1
-\nabla_{y_2}\cdot\int_{-1}^{1}\!\!\!D^\textrm{m}\nabla_{y_2}u_0^\textrm{m}
 \textrm{d}z_1
-\nabla_{y_2}\cdot\int_{-1}^{1}\!\!\!D^\textrm{m}\nabla _{z_1}u_0^\textrm{m}
 \textrm{d}z_1
\\
&
-\int_{-1}^{1}\frac{\partial}{\partial y_2}(D^\textrm{m}_{21}g(u_0^\textrm{m}))
 \textrm{d}z_1
-D^\textrm{r}\nabla U^\textrm{r}\cdot n\big|_{z_1=+1}
-D^\textrm{l}\nabla U^\textrm{l}\cdot n\big|_{z_1=-1}
=\int_{-1}^{1}F^\textrm{m}\textrm{d}z_1
.
\end{split}
\end{displaymath}
Now we integrate with respect to $y_2$ and we obtain
\begin{displaymath}
\begin{split}
\int_{0}^{1}
\int_{-1}^{1}
&
\frac{\partial u_0^\textrm{m}}{\partial T}
\textrm{d}y_2
\textrm{d}z_1
-
\int_{-1}^{1}
\left[D^\textrm{m}_{22}\frac{\partial u^\textrm{m}_0}{\partial y_2}
\right]_{y_2=0}^{y_2=1}
\textrm{d}z_1
-
\int_{-1}^{1}
\left[D^\textrm{m}_{21}\frac{\partial u^\textrm{m}_0}{\partial z_1}
\right]_{y_2=0}^{y_2=1}
\textrm{d}z_1
\\
&
-\int_{-1}^{1}\left[D^\textrm{m}_{21}g(u_0^\textrm{m})\right]_{y_2=0}^{y_2=1}
 \textrm{d}z_1
-
\int_{0}^{1}
\int_{-1}^{1}
\frac{\partial}{\partial y_2}(D^\textrm{m}_{21}g(u_0^\textrm{m}))
 \textrm{d}y_2
 \textrm{d}z_1
\\
&
-
\int_{0}^{1}
\Big[
D^\textrm{r}\nabla U^\textrm{r}\cdot n\big|_{z_1=+1}
+D^\textrm{l}\nabla U^\textrm{l}\cdot n\big|_{z_1=-1}
\Big]
 \textrm{d}y_2
=
\int_{0}^{1}
\int_{-1}^{1}
F^\textrm{m}
 \textrm{d}y_2
\textrm{d}z_1.
\end{split}
\end{displaymath}
Now, we note that the second equation in \eqref{inf070}
yields 
\begin{displaymath}
D^\textrm{m}_{21}\frac{\partial u^\textrm{m}_0}{\partial z_1}
+
D^\textrm{m}_{22}\frac{\partial u^\textrm{m}_0}{\partial y_2}
+
D^\textrm{m}_{21}g(u^\textrm{m}_0)
=0
\end{displaymath}
on $\Gamma_\textrm{h}\cap\Omega_\textrm{m}$. Recalling that $D^\textrm{m}$ and 
$u^\textrm{m}_0$ are $y_2$--periodic functions, we find the 
aforementioned jump condition 
\begin{displaymath}
\int_{0}^{1}
\Big[
D^\textrm{r}\nabla U^\textrm{r}\cdot n\big|_{z_1=+1}
+D^\textrm{l}\nabla U^\textrm{l}\cdot n\big|_{z_1=-1}
\Big]
 \textrm{d}y_2
=
\!
\int_{0}^{1}
\!\!
\int_{-1}^{1}
\Bigg[
\frac{\partial u_0^\textrm{m}}{\partial T}
-
\frac{\partial}{\partial y_2}(D^\textrm{m}_{21}g(u_0^\textrm{m}))
-
F^\textrm{m}
\Bigg]
 \textrm{d}y_2
\textrm{d}z_1.
\end{displaymath}

The relations  (\ref{inf270}) and (\ref{inf275}) provide direct access to the jump in the flux of matter when crossing the membrane. Interestingly from a modeling point of view, we can also obtain a quantitative description of the jump in concentrations  across the reduced membrane, say $\Gamma$; the situation is somehow similar to the case described in Theorem 2.4 in \cite{membrane};



\section{Numerical illustration of the finite-thickness upscaled membrane}\label{numerics}

We numerically illustrate the behavior of the finite-thickness upscaled membrane derived in Section \ref{s:duescale}. To fix a scenario, we imagine diffusion and drift of a mass concentration of gaseous CO$_2$ supposed to cross a membrane with finite thickness.

Experimental values of CO$_2$ in cells have been estimated at $d = \SI{3.5}{\centi\meter\squared\per\second}$ (cf. \cite{Kropp}). We choose diffusion coefficients around this value, i.e. $d_{1} = \SI{10}{\centi\meter\squared\per\second}$ and $d_{2} = \SI{1}{\centi\meter\squared\per\second}$, letting horizontal diffusion dominate the process.
We choose the non-linear transport term from \eqref{int000} with $b=2$.
Initially, there is no mass present, i.e. $u(t=0) = 0$.
We fix the inflow of the left boundary by choosing $u_r = \SI{5.8e-5}{\gram\per\centi\meter\cubed}$ according to \cite{Kropp} and let $u_l = 0$.
The geometry has the following dimensions: $l=\SI{1}{\centi\meter}$, $h=\SI{0.4}{\centi\meter}$, $w =\SI{0.25}{\centi\meter}$.

As $D_{22}$ lies in $L^\infty(\Omega_{\rm m}\text{-}\Omega_{\rm o})$, solving the parameter-dependent ODEs 
 \begin{equation}\label{e1}
 \frac{\partial}{\partial Y_2}\left( D_{22}(Y_1,Y_2)\frac{\partial w_1}{\partial Y_2}\right)=0
 \end{equation}
 and 
 \begin{equation}\label{e2}
 \frac{\partial}{\partial Y_2}\left[ D_{22}(Y_1,Y_2)\left(1+\frac{\partial w_2}{\partial Y_2}\right)\right]=0,
 \end{equation} 
is rather delicate since it involves distributions localized along $\partial \Omega_{\rm o}$. To handle this issue, one needs a convenient regularization of the "contrast jump". It is worth also noting that, based on \eqref{e1}-\eqref{e2}, the coefficient $D_{11}$  plays no role in the construction of the cell functions.
 Instead of smoothing the contrast, we suggest the following regularization:  Take  $\delta = \mathcal{O}(\eta)$. Find ($w_1,w_2$) such that
  \begin{eqnarray}
      \delta \frac{\partial}{\partial Y_1}(D_{11}(Y_1,Y_2)\frac{\partial}{\partial Y_1}w_1)+\frac{\partial}{\partial Y_2}(D_{22}(Y_1,Y_2)\frac{\partial}{\partial Y_2} w_1) & = & -\sqrt{\delta} \frac{\partial}{\partial Y_1}D_{11}(Y_1,Y_2), \label{pbsing1}\\
  \delta \frac{\partial}{\partial Y_1}(D_{11}(Y_1,Y_2)\frac{\partial}{\partial Y_1}w_2)+\frac{\partial}{\partial Y_2}(D_{22}(Y_1,Y_2)\frac{\partial}{\partial Y_2} w_2) & = & -\frac{\partial}{\partial Y_2}D_{22}(Y_1,Y_2).\label{pbsing2}
  \end{eqnarray}
These formulations are obtained based on (\ref{cell1}) by interpreting $\nabla_{Y_2}$ as $ \left(\begin{array}{c}  \sqrt{\delta}\frac{\partial }{\partial y_1} \\ \frac{\partial }{\partial y_2}\end{array}\right)$ instead of $\nabla_{Y_2}=\left(\begin{array}{c} 0 \\ \frac{\partial }{\partial y_2}\end{array}\right)$.
The boundary conditions needed to complete the regularized problem are described in (\ref{cell2}). This procedure appears to work well for symmetric obstacles. Note that both problems (\ref{pbsing1}) and (\ref{pbsing2}) are singular perturbations of linear elliptic PDEs. The convergence $\delta\to 0$ can be made rigorous in terms of weak solutions via a weak convergence procedure using symmetry restrictions and dimension reduction arguments.

\begin{figure}[h!]
  \includegraphics[width=0.48\textwidth]{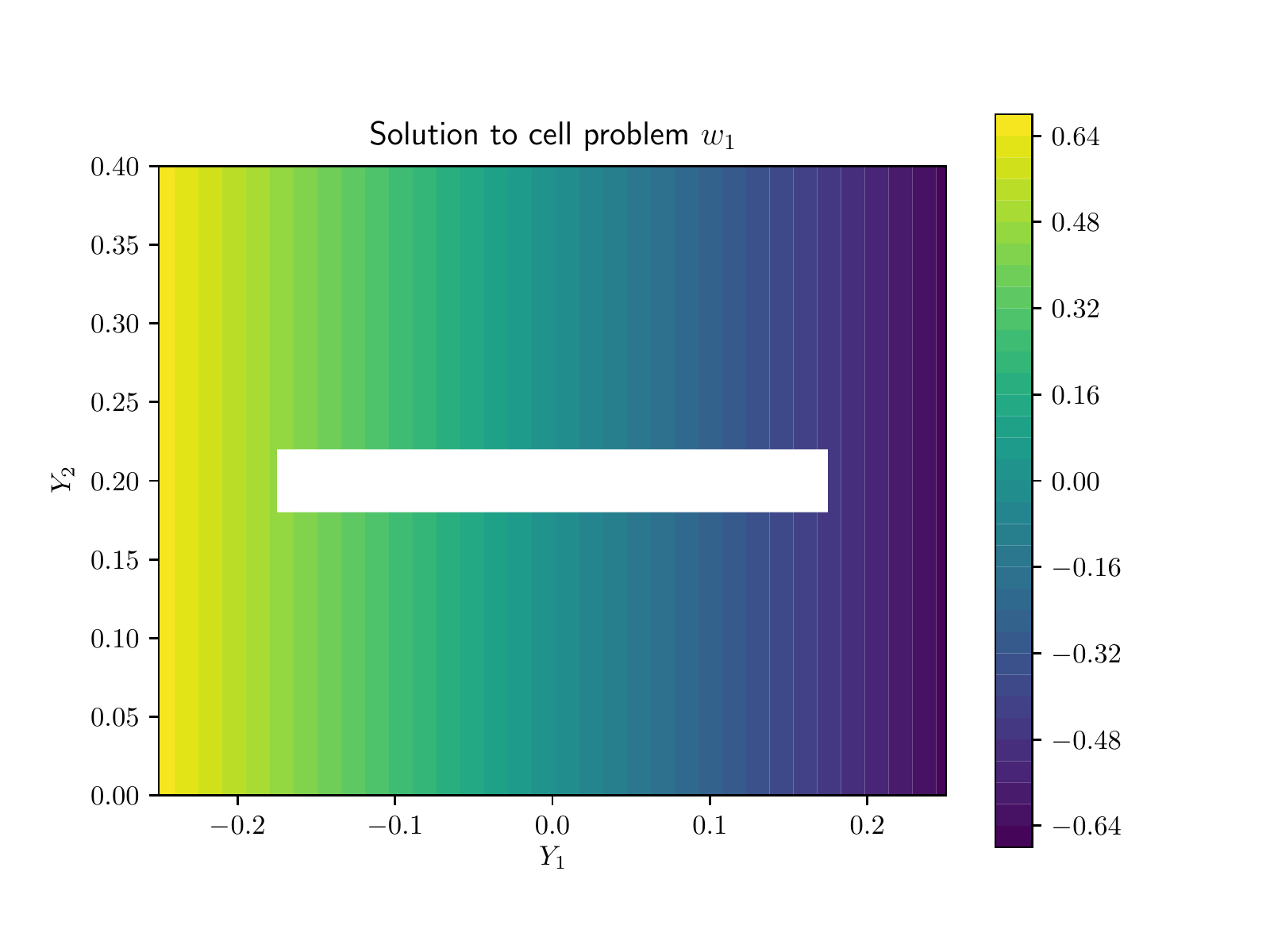}
  \includegraphics[width=.48\textwidth]{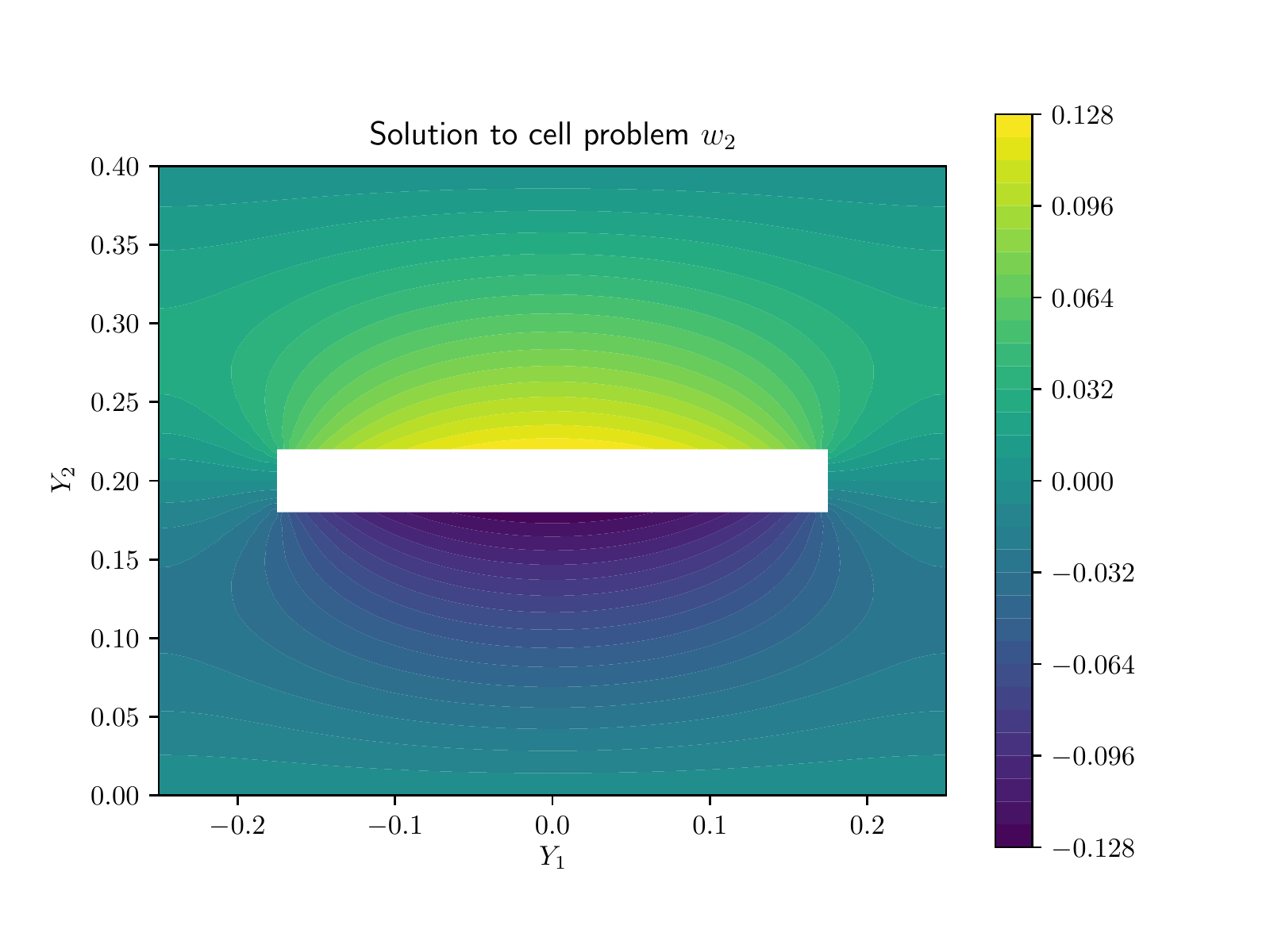}
\caption{Cell functions profiles: $w_1$ (left) and $w_2$ (right).}\label{cell-pic}
\end{figure}

To solve the cell problems  (\ref{pbsing1}) and (\ref{pbsing2}) (with corresponding boundary conditions),  we use a FEM scheme implemented in FEniCS\footnote{This is an open source platform  \texttt{FEniCS} \cite{Logg}; see \texttt{https://fenicsproject.org}.}. 
The cell problem and macroscopic equations are solved on a triangular mesh with quadratic basis functions.
We illustrate the behavior of the cell functions in Figure \ref{cell-pic}.
 
\begin{figure}[h!]
  \includegraphics[width=.50\textwidth]{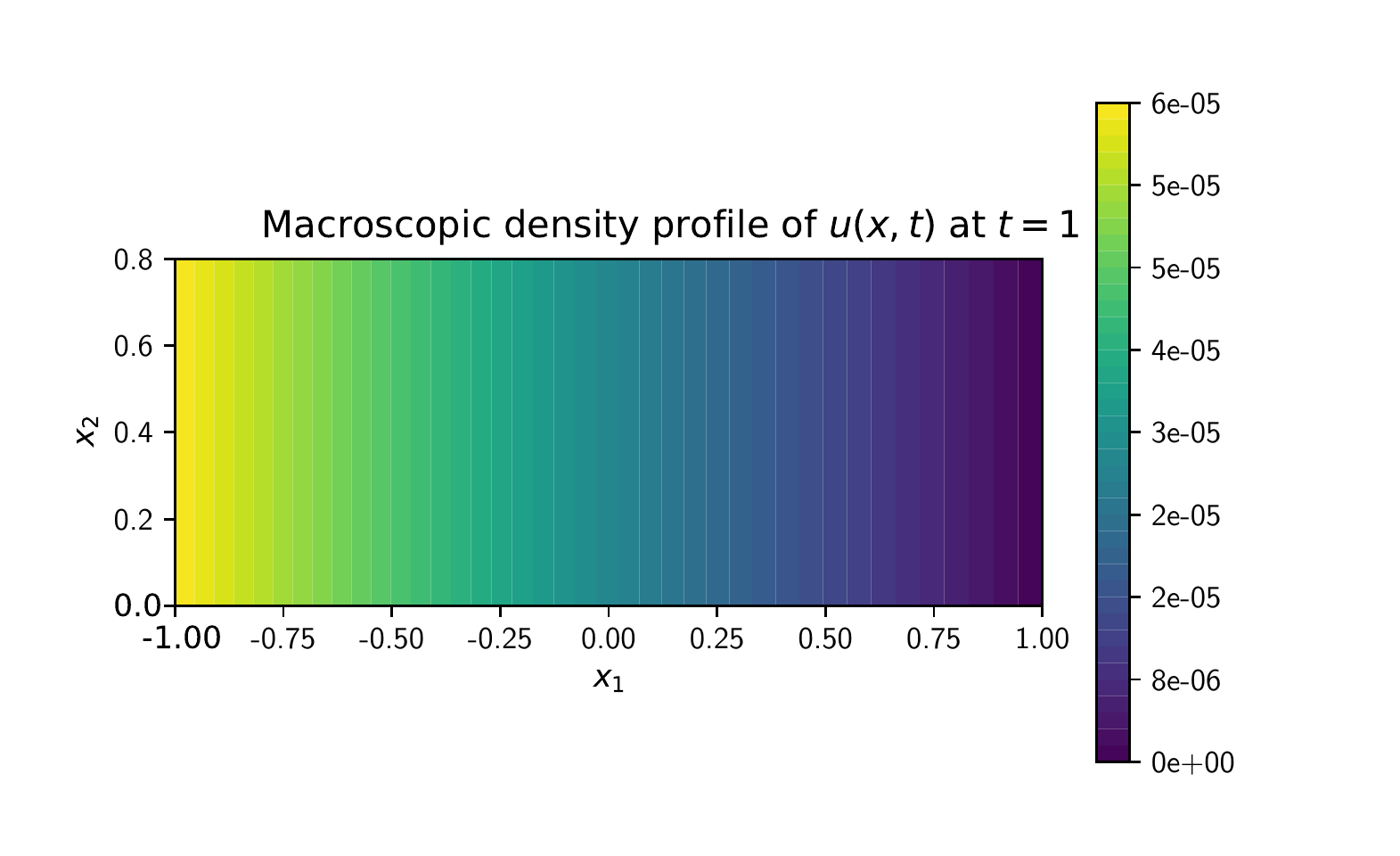}  
  \includegraphics[width=.50\textwidth]{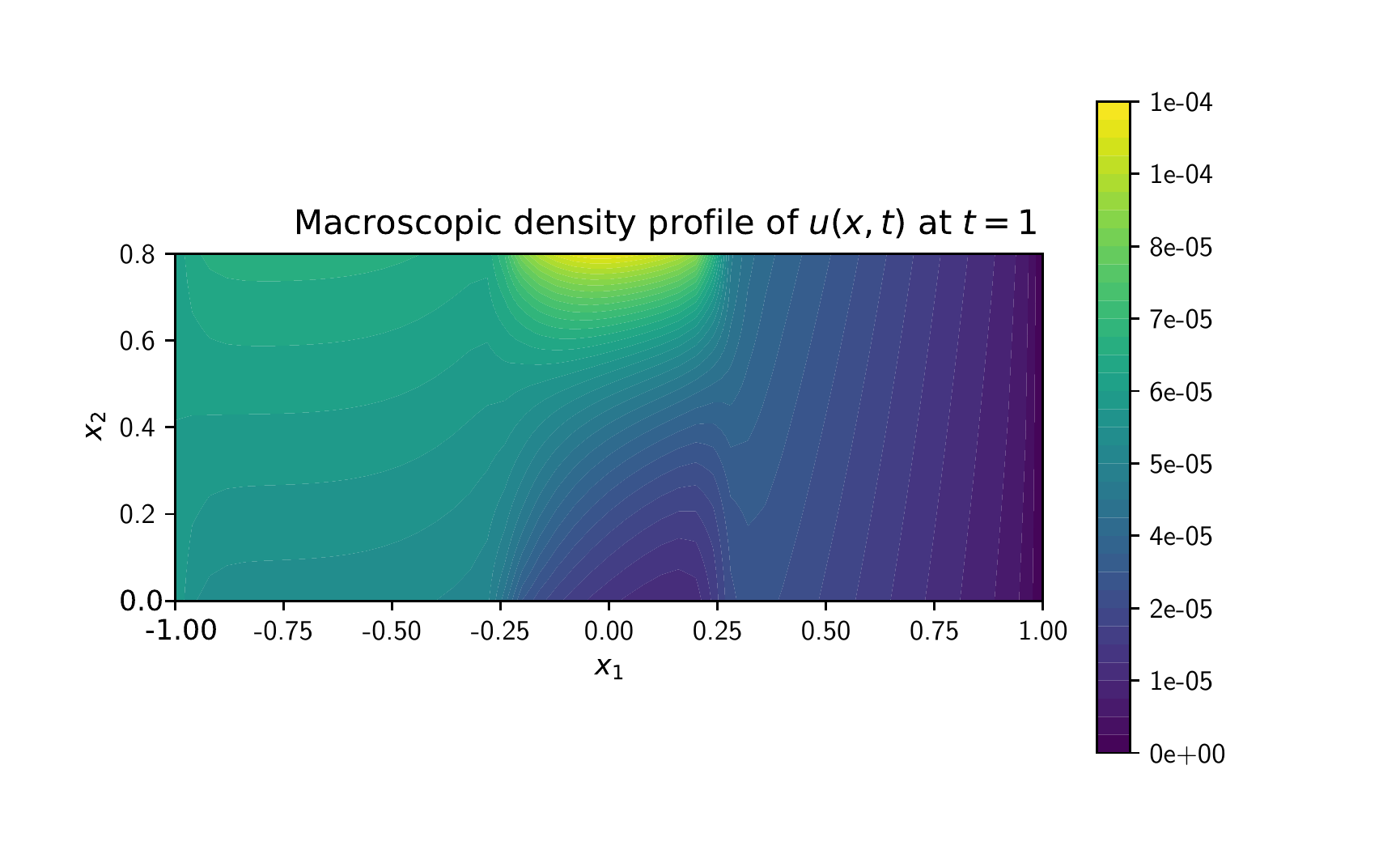}
\caption{Typical macroscopic diffusion profiles. Left:  A moderate permeability regime; Right:  Increased barrier regime exhibiting a nearly empty membrane. Interestingly, the membrane starts to behave like a barrier only in the high drift regime (i.e. for large $b$).}\label{rez0}
\end{figure}

The explicit appearance of the variable $Y_1$ in  \eqref{final01}--\eqref{final0f} needs to be removed by integrating the system  of equations with respect to the $Y_1$ variable. 
Using the transmission conditions at $\Gamma_\textrm{l}$ and $\Gamma_\textrm{r}$, the information in $\Omega_\textrm{m}$ is now linked (in a well-posed way) with equation \eqref{adim010} posed in $\Omega_\textrm{l}$ and $\Omega_\textrm{r}$, respectively.  The numerical approximations of the cell functions can now be used to compute the effective diffusion tensor
\begin{equation}
\mathbb{D}^\star:=\left(\begin{array}{cc}D_{11}^\star& D_{12}^\star \\ D_{21}^\star&D_{22}^\star\end{array}\right)=D\left(\mathbb{I}+\left(\begin{array}{cc} 0 & 0 \\ \frac{\partial w_1}{\partial Y_2}& \frac{\partial w_2}{\partial Y_2}\end{array} \right )\right),
 \end{equation}
 and hence, FEM approximations of  the upscaled diffusion-drift equation can be reached.  Note that $D^{-1}\mathbb{D}^\star$ is the so-called membrane tortuosity tensor. Typical macroscopic concentration profiles are shown in Figure \ref{rez0}. For the  chosen parameter regime, one can see that the membrane is usually permeable. Interestingly,  the efficiency of the transport through the membrane reduces when increasing the strength of the drift $b$.  Figure \ref{rez0} (right) is obtained via turning the diagonal matrix  $\mathbb{D}^\star$ into a full matrix by adding diffusion correlations. The off-diagonal entries are small $D_{12}\star = -0.05$ and $D_{21}^\star = +0.05$. Combined with a polynomial drift (of type $bu(1-u)$ with $b=54$) this causes some sort of anisotropic clogging.
 
\begin{figure}[h!]
  \includegraphics[width=.48\textwidth]{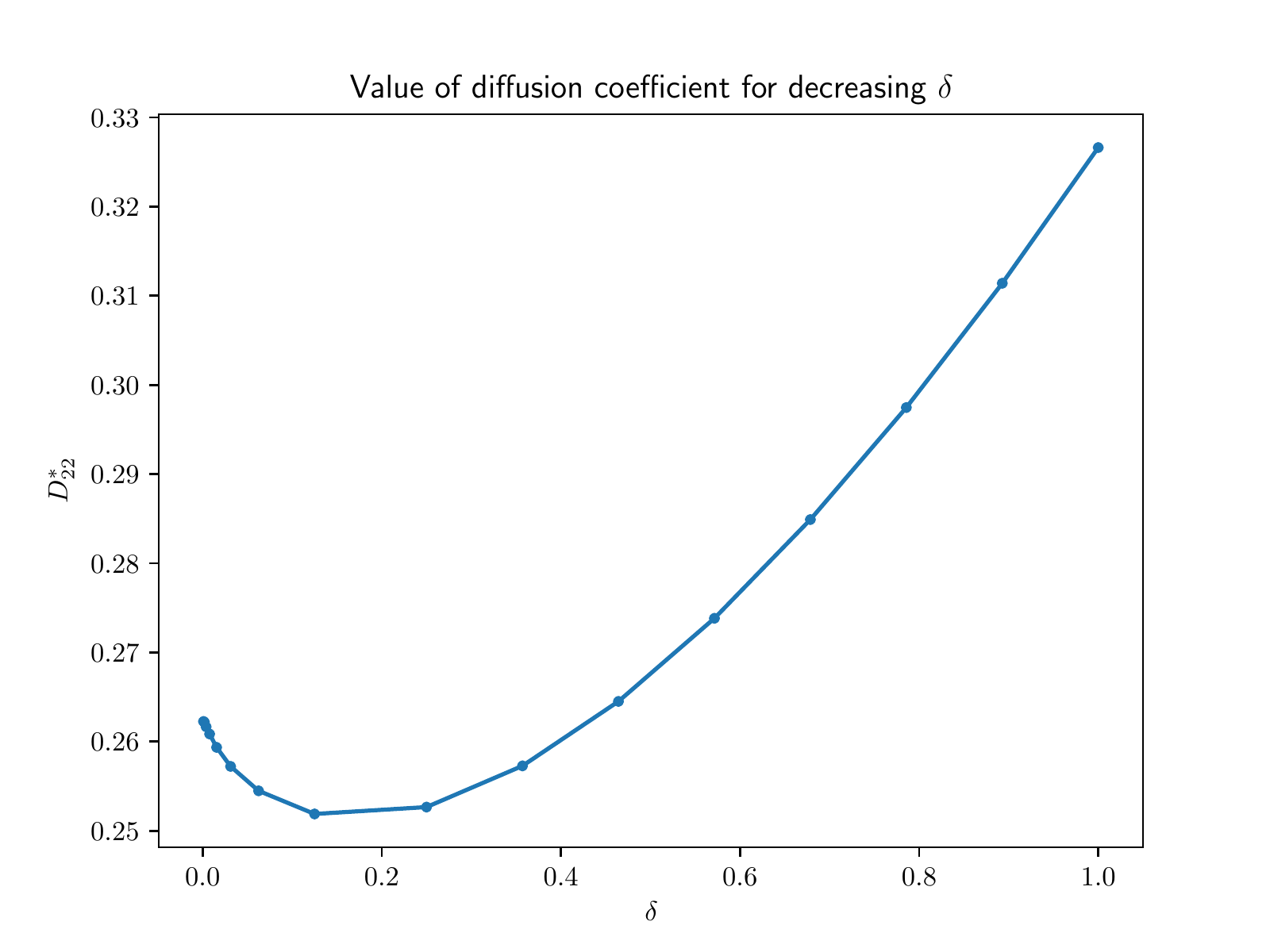}
  \includegraphics[width=0.48\textwidth]{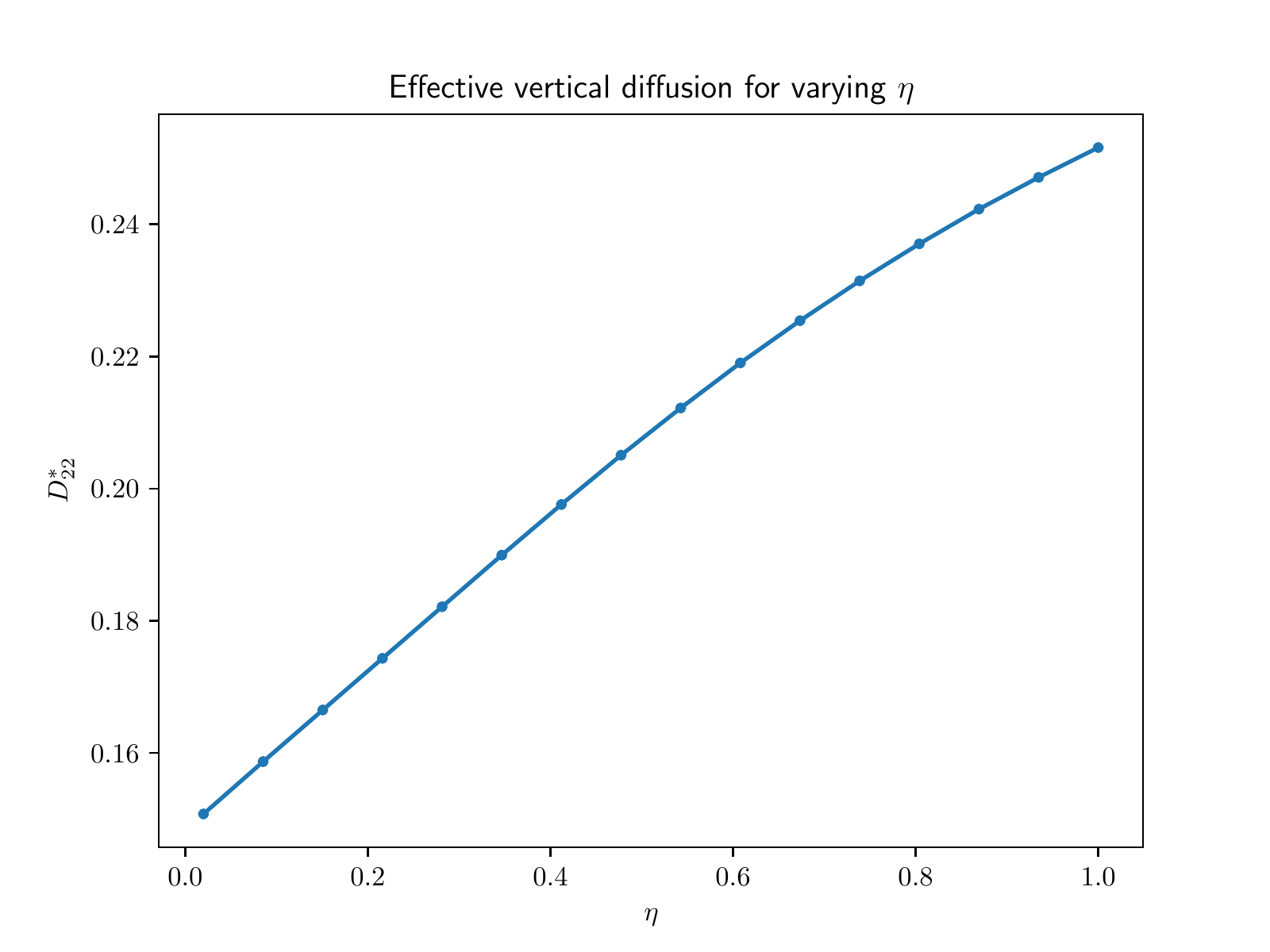}

\caption{Non-monotonicity of  $D_{22}^\star$ with respect to $\delta$, as arising in (\ref{pbsing1})--(\ref{pbsing2}). Stability of $D_{22}^\star$ with respect to the height of the periodic cell $\eta$.}\label{rez}
\end{figure}

Although the finite-thickness membrane scaling is rather standard (in the sense that the structure of the upscaled coefficients was foreseeable),  Figure \ref{rez} (left) points out an outstanding opportunity: The numerical example shows that changing the aspect ratio of the rectangular obstacle can be used as tool to optimize the membrane performance (in the spirit of shape optimization). This inspired the following key question: Is such non-monotonic behavior specific to the choice of rectangles as microstructures, or is it actually {\em generic}?

To answer this question, intensive simulations involving a large variety of shapes of microstructures need to be performed, but this is a problem by itself. For instance, the possibility of "concentration trapping" needs to be studied by, for instance, carefully considering the effect of the curvature of the micro-boundaries on the macroscopic outflux.  We will address this issue somewhere else. At this moment, relying on the stability with respect to changes in $\eta$ shown in  Figure \ref{rez} (right), we only speculate that the answer to the question is affirmative. If this were true, then, somewhat similarly to the work done in \cite{Filtration}, one can start thinking of optimizing filtration processes by searching for best-suitable microstructure shapes. This would be a useful tool for a number of engineering applications.  What the finite membrane scaling  is concerned, the optimization problem is straightforward, since it can be linked exclusively to the structure of the cell problem.  For the second scaling instead, i.e. for the infinitely-thin upscaled membrane  model, the optimization problem is not easily accessible. Here, any route towards optimizing filtration needs to take into account the structure of the limit two-scale model with nonlinear transmission condition; see \eqref{bi}--\eqref{bf}.  

\section{Discussion}\label{Discussion}

Starting from a mean-field limit of a totally asymmetric simple exclusion process (TASEP), we have investigated the problem of diffusion and non-linear drift through a composite membrane  in two specific scaling regimes. We have obtained upscaled model equations for the finite-length membrane as well as for the infinitely-thin membrane.  We can  explicitly see how the membrane microstructure affects the resulting upscaled equations and  the entries of the tensorial effective transport coefficients and our simulations show that these effects are visible at the macroscale.  From the perspective of material design, we have seen that at least what concerns the penetration of {\rm CO}$_2$ through paper, there are parameter options that can be used to optimize the membrane performance by carefully exploring the effect of the choice of the  microstructure shapes on the effective transport fluxes. 

To gain additional confidence in the model equations further investigations are needed. Two directions are more prominent: 

(i) The upscaling needs to be made mathematically rigorous. We foresee that  the two-scale convergence and boundary layer working techniques from   \cite{membrane} can be adapted to our scenario, provided one can handle the passage to the homogenization limit in the non-linear drift terms in both scalings. Additionally, the knowledge of the asymptotic expansions behind the singular perturbation (dimension reduction)--homogenization procedure can potentially be  used to derive convergence rates for the involved limiting processes.

(ii) The stochastic particle simulations from \cite{CKMSpre2016} need to be extended from the one-barrier-case to the thin composite case. Then the stationary concentration profiles and the particles residence time can be compared with findings based on the  finite element approximations of the upscaled model (both single and two-scale).  We have chosen to include solid rectangles as microstructures precisely so that the comparison between the lattice model and the upscaled evolution equations becomes possible.
 Such comparison would shed light not only on transport matters through thin porous layers (like gaseous {\rm O}$_2$ and {\rm CO}$_2$ through paper), but would also bring understanding on the effect the knowledge of the heterogeneous environments has on the stochastic dynamics of active particles (agents). 

 \paragraph{Acknowledgments.} AM and ENMC thank Prof. Rutger van Santen (Eindhoven) for fruitful discussions that have initiated this investigation. AM acknowledges a partial financial support from NWO-MPE  "Theoretical estimates of heat losses in geothermal wells" (grant No.657.014.004).  ENMC thanks FFABR 2017 financial support.


\end{document}